\documentclass[preprint,12pt]{elsarticle}

\usepackage{amsfonts}
\usepackage{graphicx}
\usepackage{epsfig}
\usepackage{geometry}
\usepackage{tikz}
\usetikzlibrary{trees,positioning,shapes}
\usepackage{subcaption}
\usepackage{amsmath}
\usepackage{amsthm}

\usepackage{xcolor, colortbl}

\newtheorem{remark}{Remark}

\newcommand{\R}{\mathbb R}

\begin{document}

\begin{frontmatter}

\title{Isogeometric Residual Minimization Method (iGRM)\\ with Direction Splitting for Non-Stationary Advection-Diffusion Problems}

\author{M. \L{}o\'{s} $^{{\textrm{(1)}}}$, 
J. Mu\~{n}oz-Matute $^{{\textrm{(2)}}}$, 
I. Muga$^{{\textrm{(3)}}}$, 
M. Paszy\'{n}ski $^{{\textrm{(1)}}}$}

\address{$^{\textrm{(1)}}$ Department of Computer Science, \\ AGH University of Science and Technology,
Krakow, Poland \\
e-mail: paszynsk@agh.edu.pl \\
e-mail: marcin.los.91@gmail.com}

\address{$^{\textrm{(2)}}$ The University of the Basque Country, \\ Bilbao, Spain \\
e-mail: judith.munozmatute@gmail.com}

\address{$^{\textrm{(3)}}$ Pontificia Universidad Cat\'{o}lica Valpara\'iso, Chile \\
e-mail: ignacio.muga@pucv.cl}

\begin{abstract}
In this paper, we propose a novel computational implicit method, which we call Isogeometric Residual Minimization (iGRM) with direction splitting. The method mixes the benefits resulting from isogeometric analysis, implicit dynamics, residual minimization, and alternating direction solver. We utilize tensor product B-spline basis functions in space, implicit second order time integration schemes, residual minimization in every time step, and we exploit Kronecker product structure of the matrix to employ linear computational cost alternating direction solver. We implement an implicit time integration scheme and apply, for each space-direction, a stabilized mixed method based on residual minimization. We show that the resulting system of linear equations has a Kronecker product structure, which results in a linear computational cost of the direct solver,
even using implicit time integration schemes together with the stabilized mixed formulation. 
We test our method on three advection-diffusion computational examples, including model ``membrane'' problem, the circular wind problem, and the simulations modeling pollution propagating from a chimney.
\end{abstract}
	

\end{frontmatter}

\section{Introduction}

From one side, the alternating directions method (ADS) was introduced in \cite{ADS1,ADS2,ADS3,ADS4} to solve parabolic and hyperbolic problems using finite differences. A modern version of this method solves different classes of problems \cite{Minev1,Minev2}.

From the other side, Isogeometric Analysis (IGA) \cite{IGA} bridges the gap between the Computer Aided Design (CAD) and Computer Aided Engineering (CAE) communities. The idea of IGA is to apply spline basis functions \cite{NURBS} for finite element method (FEM) simulations. The ultimate goal is to perform engineering analysis directly to CAD models without expensive remeshing and recomputations.
IGA has multiple applications in time-dependent simulations, including phase-field models \cite{Dede2010,Dede2011}, phase-separation simulations with application to cancer growth simulations \cite{Gomez2008,Gomez2010}, wind turbine aerodynamics \cite{Hsu:2011}, incompressible hyper-elasticity \cite{Duddu:2011}, turbulent flow simulations \cite{Chang:2011}, transport of drugs in cardiovascular applications \cite{Hossain:2011}, or the blood flow simulations and drug transport in arteries simulations \cite{Bazilevs2006,Bazilevs2007,Calo:2008}.

The combination of ADS methods and IGA was already applied \cite{CG1,CG2,Gao2014} for fast solution of the $L^2$ projection problem with isogeometric FEM. Indeed
an explicit time integration scheme with isogeometric discretization is equivalent to the solution of a sequence of isogeometric $L^2$ projections, and the direction splitting method is a critical tool for a rapid computation of the aforementioned projections. This idea was successfully used for fast explicit dynamics simulations \cite{ED1,ED2,ED3,ED4,ED5}. 


The stability of a numerical method based on Petrov-Galerkin discretizations of a general weak form relies on the famous discrete inf-sup condition (see, e.g.,~\cite{ErnGuermondBOOK2004}):
\begin{equation}
 \sup_{v_h \in V_h} {{|b(u_h,v_h)|} \over {\|v_h\|_{V_h}}} \geq \gamma_h \| u_h \|_{U_h}, \quad\forall u_h\in U_h \subset U,
 \label{eq9}
\end{equation}
where $b(\cdot,\cdot)$ is a given bilinear form, $\gamma_h>0$ is the discrete inf-sup (stability) constant, and $(U_h,V_h)$ are the chosen discrete trial and test spaces, respectively. However, to find a {\it compatible pair} $(U_h,V_h)$ satisfying \eqref{eq9} with nice stability constant is far from being a simple task, especially when dealing with isogeometric basis functions. One way to somehow circumvent the last requirement is to perform a {\it Residual Minimization} method, that instead, aims to compute $u_h\in U_h$ such that:
\begin{equation}\label{eq:res_min}
u_h=\arg\!\!\!\min_{w_h\in U_h}\|b(w_h,\cdot)-l(\cdot)\|_{V'_h}\,\,,
\end{equation}
where $l(\cdot)$ is a given linear functional that corresponds to the right-hand side of the weak formulation. Indeed, if the bilinear form $b(\cdot,\cdot)$ is inf-sup stable at the continuous (infinite dimensional) level, then there is always a {\it hope} that, by smartly increasing the dimension of the discrete test space, discrete stability is reached at some point.

There is consistent literature on residual minimization methods, especially for convection-diffusion problems \cite{evans,dahmen,stevenson}, where it is well known that the lack of stability is the main issue to overcome. In particular, the class of DPG methods \cite{dh1,dh2} aim to obtain a practical approach to solve the mixed system derived from~\eqref{eq:res_min} by breaking the test spaces (at the expense of introducing a hybrid formulation). Breaking the test spaces is essential for a practical application of the Schur complement technique and the consequent reduction problem size. In this work, we have taken a different path: instead of using hybrid formulations and broken spaces, we take advantage of the tensor product structure of the computational mesh and Kronecker product structure of matrices, to obtain linear computational cost ${\cal O}(N)$ solver.


There are many methods constructing test functions providing better stability of the method for given class of problems \cite{t4,t5,t6,t7}.
In 2010 the Discontinuous Petrov-Galerkin (DPG) method was proposed, with the modern summary of the method described in \cite{t8,t9}.
The key idea of the DPG method is to construct the optimal test functions ``on the fly'', element by element.
The DPG automatically guarantee the numerical stability of difficult computational problems,
thanks to the automatic selection of the optimal basis functions.
The DPG method is equivallent to the residual minimization method \cite{t8}.
The DPG is a practical way to implement the residual minimization method, when the computational cost of the global solution is expensive (non-linear).

In this paper, we seek to develop a new computational paradigm to obtain stable and accurate solutions of time-dependent problems, using  implicit schemes, which we call isogeometric Residual Minimization (iGRM), with the following unique features:
\begin{enumerate}
\item On tensor product meshes, linear computational cost O(N) of the direct solver solution,
\item Unconditional stability of the implicit time integration scheme,
\item Unconditional stability in the spatial domain.
\end{enumerate}
This new paradigm may have a significant impact on the way how the computational mechanics community performs simulations of time-dependent problems on tensor product grids. 
The methodology we propose brings together the benefits of the Discontinuous Petrov-Galerkin method (DPG), Isogeometric Analysis (IGA), and Alternating Direction Implicit solvers (ADI).

One of the advantages of the isogeometric residual minimization is a possibility of preserving arbitrary smoothness of the approximation, which is difficult to obtain with alternative DPG method.
The main point of using residual minimization methods in this work, is to derive a stabilization method that does not break the Kronecker product structure. 
In particular, we increase the test spaces by increasing the order of the B-splines, and/or decreasing the continuity.

One of the motivations for the broken spaces used in the Discontinuous Petrov-Galerkin (DPG) method is to enable a local static condensation. 
We do not perform the static condensation in our method, we rather solve the global problem once, exploiting the Kronecker product structure of the entire system, including the Gramm matrix resulting from the inner product at the test space, and the weak form matrix. 
The Kronecker product structure of the entire system is preserved even after the residual minimization method is applied. Each one-dimensional block of this system has a form presented in Figure 1. For this kind of matrix, we can factorize it with a linear computational cost ${\cal O}(N)$.
The application of the broken spaces requires the introduction of the fluxes and traces on the boundaries of elements. These terms may be problematic to preserve the Kronecker product structure of the entire system.
The hybridization process introduces new variables and equations that may break the Kronecker product structure. 
In our case, we have a linear computational cost solver for the residual minimization method without hybridization resulting in additional traces and/or fluxes variables.

We approach the residual minimization using its saddle point (mixed) formulation, e.g., as described in~\cite{evans}.   
Several discretization techniques are particular incarnations of this wide-class of residual minimization methods. These include: the \emph{least-squares finite element method}~\cite{LSFEM}, the \emph{discontinuous Petrov-Galerkin method (DPG) with optimal test functions}~\cite{t8}, or the \emph{variational stabilization method}~\cite{VSM}. 

We apply our method for the solution of the advection-diffusion problems.

For our method to work, we need to preserve the Kronecker product structure of the matrix, and this includes the Jacobian from changing of the variables from the physical geometry into the regular patch of elements. In general, it is possible to deal with rectangular, spherical or cylindrical domains. For other geometries of the patch of elements, we could still use our fast solver as a linear cost preconditioner for an iterative solver, just by replacing the jacobian by a Kronecker product approximation, following the ideas of \cite{CG2}.

There are several possibilities \cite{Norm1,dh1,dh2} to choose the inner product for the advection-diffusion problems, see e.g. section 2.2 in \cite{dh2}. These norms, however, usually include $L^2$ norm of $u$, $L^2$ norm of $\nabla u$, $L^2$ norm of $\beta\cdot \nabla u$, etc.
In our method, we do not use hybrid formulation, because our target is not to obtain a block diagonal matrix but a system that can be solved in a linear cost taking advantage of the Kronecker product structure. Thus, in our formulation, we do not include traces and/or fluxes, and we do not break the test spaces. In each time step, we solve the differential problem that results from splitting of the whole operator into two sub-operators, each of them having only derivatives in one direction. We perform similar splitting for the inner product term. This is critical to preserve the Kronecker product structure and thus the linear cost of the solver.
So there is not much choice in the selection of the norms in our case. 
It is possible to try the weighted norms resulting in the inner product split into components like $(u,u)+\alpha (\frac{\partial u}{\partial x_i},\frac{\partial u}{\partial x_i})$, where $\alpha$ is related to diffusion and/or advection terms.

The Petrov-Galerkin type of stabilization was applied to simulations of time-dependent problems for the first order transport equations \cite{DahmenDPG}.
The main original contribution of our work is to provide a linear ${\cal O}(N)$ computational cost solver for the factorization of the stabilized system in every time step.
As mentioned for the least-square methods \cite{DGmethod}, the structure of the Gramm matrix and the factorization cost depends on the test spaces used. If the spaces are broken, we get the block-diagonal matrix, and the Schur complements can be computed locally over each element. If the spaces are continuous, we have to solve the global system. In our case, the setup is similar, and we aim to provide the Kronecker product structure of the global system without breaking the test spaces, to have the linear cost factorization.
Smart selections of the test norm can impact the convergence of the residual minimization method \cite{Alternative}. Additionally, different weak or ultra-weak formulations have an impact on the convergence as well as \cite{RMmixed}.
Different time integration schemes have been also considered in the context of DPG method for heat transfer equations \cite{HeatDPG}.
In our work, we focus on inner products, formulations, and time integration schemes, that may lead to the Kronecker product structure of the matrix, the linear computational cost, and the second-order in time accuracy.


\section{Model problem}\label{S2}

Let $\Omega=\Omega_{x}\times\Omega_{y}\subset \R^{2}$ a bounded domain and $I=(0,T]\subset\R$, we consider the two-dimensional \textit{linear advection-diffusion equation}

\begin{equation}\label{strong}
\displaystyle{ \left\{
\begin{split}
u_{t}-\nabla \cdot(\alpha\nabla u)+\beta\cdot \nabla u=&\;f&\mbox{in}&\;\Omega\times I,\\
u=&\;0&\mbox{on}&\;\Gamma\times I,\\
u(0)=&\;u_{0}\;&\mbox{in}&\;\Omega,\\
\end{split}
\right.} 
\end{equation}
where $\Omega_{x}$ and $\Omega_{y}$ are intervals in $\R$. Here, $u_{t}:=\partial u/\partial t$, $\Gamma=\partial\Omega$ denotes the boundary of the spatial domain $\Omega$, $f:\Omega\times I\longrightarrow\R$ is a given source and $u_{0}:\Omega\longrightarrow\R$ is a given initial condition. 
We consider constant diffusivity $\alpha$ and a velocity field $\beta=[\beta_{x}\;\;\beta_{y}].$

\section{Operator splitting}\label{S3}
We split the advection-diffusion operator $\mathcal{L}u=-\nabla \cdot(\alpha\nabla u)+\beta\cdot \nabla u$ as $\mathcal{L}u=\mathcal{L}_{1}u+\mathcal{L}_{2}u$ where
$$\mathcal{L}_{1}u:=-\alpha\frac{\partial^2 u}{\partial x^{2}}+\beta_{x}\frac{\partial u}{\partial x},\;\;\;\mathcal{L}_{2}u:=-\alpha\frac{\partial^2 u}{\partial y^{2}}+\beta_{y}\frac{\partial u}{\partial y}.$$
Based on this operator splitting, we consider different Alternating Direction Implicit (ADI) schemes to discretize problem (\ref{strong}).

First, we perform an uniform partition of the time interval $\bar{I}=[0,T]$ as
$$0=t_{0}<t_{1}<\ldots<t_{N-1}<t_{N}=T,$$
and denote $\tau:=t_{n+1}-t_{n},\;\forall n=0,\ldots,N-1$.

\section{Peaceman-Rachford scheme}\label{S4} 
In the Peaceman-Rachford scheme \cite{ADS1,ADS2}, we integrate the solution from time step $t_{n}$ to $t_{n+1}$ in two substeps as follows:
\begin{equation}\label{PR}
\displaystyle{\left\{
\begin{split}
\frac{u^{n+1/2}-u^{n}}{\tau/2}+\mathcal{L}_{1}u^{n+1/2}&=f^{n+1/2}-\mathcal{L}_{2}u^{n},\\
\frac{u^{n+1}-u^{n+1/2}}{\tau/2}+\mathcal{L}_{2}u^{n+1}&=f^{n+1/2}-\mathcal{L}_{1}u^{n+1/2}.\\
\end{split}
\right.}
\end{equation}
The variational formulation of scheme (\ref{PR}) is 
\begin{equation}\label{varPR}
\displaystyle{\left\{
\begin{split}
&(u^{n+1/2},v)+\frac{\tau}{2}\left(\alpha\frac{\partial u^{n+1/2}}{\partial x},\frac{\partial v}{\partial x}\right)+\frac{\tau}{2}\left(\beta_{x}\frac{\partial u^{n+1/2}}{\partial x},v\right)=\\
&=(u^{n},v)-\frac{\tau}{2}\left(\alpha\frac{\partial u^{n}}{\partial y},\frac{\partial v}{\partial y}\right)-\frac{\tau}{2}\left(\beta_{y}\frac{\partial u^{n}}{\partial y},v\right)+\frac{\tau}{2}(f^{n+1/2},v),\\
&\\
&(u^{n+1},v)+\frac{\tau}{2}\left(\alpha\frac{\partial u^{n+1}}{\partial y},\frac{\partial v}{\partial y}\right)+\frac{\tau}{2}\left(\beta_{y}\frac{\partial u^{n+1}}{\partial y},v\right)=\\
&=(u^{n+1/2},v)-\frac{\tau}{2}\left(\alpha\frac{\partial u^{n+1/2}}{\partial x},\frac{\partial v}{\partial x}\right)-\frac{\tau}{2}\left(\beta_{x}\frac{\partial u^{n+1/2}}{\partial x},v\right)+\frac{\tau}{2}(f^{n+1/2},v),\\
\end{split}
\right.}
\end{equation}
where $(\cdot,\cdot)$ denotes the inner product of $L^{2}(\Omega)$. Finally, expressing problem (\ref{varPR}) in matrix form we have
\begin{equation}\label{matrixPR}
\displaystyle{\left\{
\begin{split}
&\left[M^{x}+\frac{\tau}{2}(K^{x}+G^{x})\right]\otimes M^{y}u^{n+1/2}=M^{x}\otimes\left[M^{y}-\frac{\tau}{2}(K^{y}+G^{y})\right]u^{n}+\frac{\tau}{2}F^{n+1/2},\\
&M^{x}\otimes\left[M^{y}+\frac{\tau}{2}(K^{y}+G^{y})\right]u^{n+1}=\left[M^{x}-\frac{\tau}{2}(K^{x}+G^{x})\right]\otimes M^{y}u^{n+1/2}+\frac{\tau}{2}F^{n+1/2},
\end{split}
\right.}
\end{equation}
where $M^{x,y}$, $K^{x,y}$ and $G^{x,y}$ are the 1D mass, stiffness and advection matrices, respectively. 

\section{Strang splitting scheme}\label{S5} 
In the Strang splitting scheme \cite{Strang,Oden} we divide problem $u_{t}+\mathcal{L}u=f$ into two subproblems as follows:
\begin{equation}\label{subStrang}
\displaystyle{\left\{
\begin{split}
P_{1}:&\;u_{t}+\mathcal{L}_{1}u=f,\\
P_{2}:&\;u_{t}+\mathcal{L}_{2}u=0,\\
\end{split}
\right.}
\end{equation}
and the scheme integrates the solution from $t_{n}$ to $t_{n+1}$ in three substeps:
\begin{equation}\label{schemeStrang}
\displaystyle{\left\{
\begin{split}
\mbox{Solve}\;P_{1}:&\;u_{t}+\mathcal{L}_{1}u=f,\;\mbox{in}\;(t_{n},t_{n+1/2}),\\
\mbox{Solve}\;P_{2}:&\;u_{t}+\mathcal{L}_{2}u=0,\;\mbox{in}\;(t_{n},t_{n+1}),\\
\mbox{Solve}\;P_{1}:&\;u_{t}+\mathcal{L}_{1}u=f,\;\mbox{in}\;(t_{n+1/2},t_{n+1}),\\
\end{split}
\right.}
\end{equation}
and we can employ different methods in each substep of (\ref{schemeStrang}).
\subsection{Backward Euler}
If we select the Backward Euler method in (\ref{schemeStrang}), we obtain
\begin{equation}\label{BE}
\displaystyle{\left\{
\begin{split}
&\frac{u^{*}-u^{n}}{\tau/2}+\mathcal{L}_{1}u^{*}=f^{n+1/2},\\
&\frac{u^{**}-u^{*}}{\tau}+\mathcal{L}_{2}u^{**}=0,\\
&\frac{u^{n+1}-u^{**}}{\tau/2}+\mathcal{L}_{1}u^{n+1}=f^{n+1}.\\
\end{split}
\right.}
\end{equation}
where $u^{*}$ and $u^{**}$ denotes intermediate solutions as defined by the Strang integration scheme.

The variational formulation of scheme (\ref{BE}) is 
\begin{equation}\label{varBE}
\displaystyle{\left\{
\begin{split}
&(u^{*},v)+\frac{\tau}{2}\left(\alpha\frac{\partial u^{*}}{\partial x},\frac{\partial v}{\partial x}\right)+\frac{\tau}{2}\left(\beta_{x}\frac{\partial u^{*}}{\partial x},v\right)=(u^{n},v)+\frac{\tau}{2}(f^{n+1/2},v),\\
&(u^{**},v)+\tau\left(\alpha\frac{\partial u^{**}}{\partial y},\frac{\partial v}{\partial y}\right)+\tau\left(\beta_{y}\frac{\partial u^{**}}{\partial y},v\right)=(u^{*},v),\\
&(u^{n+1},v)+\frac{\tau}{2}\left(\alpha\frac{\partial u^{n+1}}{\partial x},\frac{\partial v}{\partial x}\right)+\frac{\tau}{2}\left(\beta_{x}\frac{\partial u^{n+1}}{\partial x},v\right)=(u^{**},v)+\frac{\tau}{2}(f^{n+1},v),\\
\end{split}
\right.}
\end{equation}
and scheme (\ref{varBE}) in matrix form is
\begin{equation}\label{matrixBE}
\displaystyle{\left\{
\begin{split}
&\left[M^{x}+\frac{\tau}{2}(K^{x}+G^{x})\right]\otimes M^{y}u^{*}=M^{x}\otimes M^{y}u^{n}+\frac{\tau}{2}F^{n+1/2},\\
&M^{x}\otimes\left[M^{y}+\tau(K^{y}+G^{y})\right]u^{**}=M^{x}\otimes M^{y}u^{*},\\
&\left[M^{x}+\frac{\tau}{2}(K^{x}+G^{x})\right]\otimes M^{y}u^{n+1}=M^{x}\otimes M^{y}u^{**}+\frac{\tau}{2}F^{n+1}.\\
\end{split}
\right.}
\end{equation}

\subsection{Crank-Nicolson}
If we select the Crank-Nicolson method in (\ref{schemeStrang}), we obtain
\begin{equation}\label{CN}
\displaystyle{\left\{
\begin{split}
&\frac{u^{*}-u^{n}}{\tau/2}+\frac{1}{2}(\mathcal{L}_{1}u^{*}+\mathcal{L}_{1}u^{n})=\frac{1}{2}(f^{n+1/2}+f^{n}),\\
&\frac{u^{**}-u^{*}}{\tau}+\frac{1}{2}(\mathcal{L}_{2}u^{**}+\mathcal{L}_{2}u^{*})=0,\\
&\frac{u^{n+1}-u^{**}}{\tau/2}+\frac{1}{2}(\mathcal{L}_{1}u^{n+1}+\mathcal{L}_{1}u^{**})=\frac{1}{2}(f^{n+1}+f^{n+1/2}).\\
\end{split}
\right.}
\end{equation}

The variational formulation of scheme (\ref{CN}) is 
\begin{equation}\label{varCN}
\displaystyle{\left\{
\begin{split}
&(u^{*},v)+\frac{\tau}{4}\left(\alpha\frac{\partial u^{*}}{\partial x},\frac{\partial v}{\partial x}\right)+\frac{\tau}{4}\left(\beta_{x}\frac{\partial u^{*}}{\partial x},v\right)=\\
&=(u^{n},v)-\frac{\tau}{4}\left(\alpha\frac{\partial u^{n}}{\partial x},\frac{\partial v}{\partial x}\right)-\frac{\tau}{4}\left(\beta_{x}\frac{\partial u^{n}}{\partial x},v\right)+\frac{\tau}{4}(f^{n+1/2}+f^{n},v),\\
&\\
&(u^{**},v)+\frac{\tau}{2}\left(\alpha\frac{\partial u^{**}}{\partial y},\frac{\partial v}{\partial y}\right)+\frac{\tau}{2}\left(\beta_{y}\frac{\partial u^{**}}{\partial y},v\right)=\\
&=(u^{*},v)-\frac{\tau}{2}\left(\alpha\frac{\partial u^{*}}{\partial y},\frac{\partial v}{\partial y}\right)-\frac{\tau}{2}\left(\beta_{y}\frac{\partial u^{*}}{\partial y},v\right),\\
&\\
&(u^{n+1},v)+\frac{\tau}{4}\left(\alpha\frac{\partial u^{n+1}}{\partial x},\frac{\partial v}{\partial x}\right)+\frac{\tau}{4}\left(\beta_{x}\frac{\partial u^{n+1}}{\partial x},v\right)=\\
&=(u^{**},v)-\frac{\tau}{4}\left(\alpha\frac{\partial u^{**}}{\partial x},\frac{\partial v}{\partial x}\right)-\frac{\tau}{4}\left(\beta_{x}\frac{\partial u^{**}}{\partial x},v\right)+\frac{\tau}{4}(f^{n+1}+f^{n+1/2},v).\\
\end{split}
\right.}
\end{equation}
and scheme (\ref{varCN}) in matrix form is
\begin{equation}\label{matrixCN}
\displaystyle{\left\{
\begin{split}
&\left[M^{x}+\frac{\tau}{4}(K^{x}+G^{x})\right]\otimes M^{y}u^{*}=\left[M^{x}-\frac{\tau}{4}(K^{x}+G^{x})\right]\otimes M^{y}u^{n}+\frac{\tau}{4}(F^{n+1/2}+F^{n}),\\
&M^{x}\otimes\left[M^{y}+\frac{\tau}{2}(K^{y}+G^{y})\right]u^{**}=M^{x}\otimes\left[M^{y}-\frac{\tau}{2}(K^{y}+G^{y})\right]u^{*},\\
&\left[M^{x}+\frac{\tau}{4}(K^{x}+G^{x})\right]\otimes M^{y}u^{n+1}=\left[M^{x}-\frac{\tau}{4}(K^{x}+G^{x})\right]\otimes M^{y}u^{**}+\frac{\tau}{4}(F^{n+1}+F^{n+1/2}).\\
\end{split}
\right.}
\end{equation}


\section{Isogeometric residual minimization method}

In all the above methods, in every time step we solve the following problem:
Find $u \in U$ such as
\begin{equation}
b\left(u,v\right)=l\left(v\right) \quad \forall v\in V,
\label{eq:bl}
\end{equation}
\begin{equation}
b\left(u,v\right)= \left(u,v\right) + \Delta t \left( \left(\beta_i \frac {\partial u}{\partial x_i}, v\right) + \alpha_i\left(\frac{\partial u}{\partial x_i}, \frac{\partial v}{\partial x_i}\right)\right).
\label{eq:b}
\end{equation}
where $\Delta t=\tau/2$ for the Peaceman-Rachford, $\Delta t=\tau/2$ for the  Strang method with backward Euler, and $\Delta t=\tau/4$ for the Strang method with Crank-Nicolson scheme.
The right-hand-side $l\left(w,v\right)$ depends on the selected time-integration scheme, e.g. for the Strang method with backward Euler it is
\begin{equation}
l\left(w,v\right)  = \left(w+\Delta tf,v\right) \quad \forall v\in V.
\label{eq:l}
\end{equation}

In our advection-diffusion problem we seek the solution in space
\begin{equation}
U = V =\left\{ v: \int_{\Omega} \left(v^2+\left(\frac{\partial v}{\partial x_i}\right)^2 \right) <\infty \right\}.
\label{eq:inner}
\end{equation}
where $i=1,2$ denotes the spatial directions.

The inner product in $V$ is defined as
\begin{equation}
\left(u,v\right)_V=\left(u,v\right)_{L_2}+\left(\frac{\partial u}{\partial x_{\underline{i}}},\frac{\partial v}{\partial x_{\underline{i}}}\right)_{L_2},
\end{equation}
where $i=1,2$ depending on the first or the second sub-step in the alternating directions method, and we do not use here the Einstein convention.

For a weak problem we define the operator $B: U \rightarrow V'$ such as 
$<Bu,v>_{V'\times V}=b\left(u,v\right)$.

\begin{equation}
 B: U \rightarrow V',
\end{equation}
such that
\begin{equation}
 \langle Bu , v \rangle_{V' \times V} = b(u,v),
\end{equation}
so we can reformulate the problem as
\begin{equation}
 Bu - l = 0.
\end{equation}
We wish to minimize the residual
\begin{equation}
 u_h = \mathrm{argmin}_{w_h \in U_h} {1 \over 2} \| Bw_h - l \|_{V'}^2.
\end{equation}
We introduce the Riesz operator being the isometric isomorphism
\begin{equation}
 R_V \colon V \ni v \rightarrow (v,.) \in V'.
\end{equation}
We can project the problem back to $V$
\begin{equation}
 u_h = \mathrm{argmin}_{w_h \in U_h} {1 \over 2} \| R_V^{-1} (Bw_h - l) \|_V^2.
\end{equation}
The minimum is attained at $u_h$ when the G\^ateaux derivative is equal to $0$ in all directions:
\begin{equation}
 \langle R_V^{-1} (Bu_h - l), R_V^{-1}(B\, w_h) \rangle_V = 0 \quad \forall \, w_h \in U_h.
\end{equation}
We define the error representation function $r=R_V^{-1}(Bu_h-l)$ and our problem is reduced to 
\begin{equation}
 \langle r , R_V^{-1} (B\,  w_h ) \rangle = 0 \quad \forall \, w_h \in U_h,
\end{equation}
which is equivalent to 
\begin{equation}
 \langle Bw_h, r  \rangle = 0 \quad \quad \forall w_h \in U_h.
\end{equation}
From the definition of the residual we have 
\begin{equation}
(r,v)_V=\langle B u_h-l,v \rangle \quad \forall v\in V.
\end{equation}
Our problem reduces to the following semi-infinite problem: Find $(r,u_h)_{V\times U_h}$ such as 
\begin{eqnarray}
\begin{aligned}
(r,v)_V - \langle B u_h-l,v \rangle &= 0 \quad  \forall v\in V \\
\langle Bw_h,r\rangle &= 0 \quad  \forall w_h \in U_h.
\end{aligned}
\end{eqnarray}
We discretize the test space $V_m \in V$ to get the discrete problem: 
Find $(r_m,u_h)_{V_m\times U_h}$ such as 
\begin{eqnarray}
\begin{aligned}
(r_m,v_m)_{V_m} - \langle B u_h-l,v_m \rangle &= 0 \quad  \forall v\in V_m \\
\langle Bw_h,r_m\rangle &= 0 \quad \forall w_h \in U_h.
\end{aligned}
\label{eq:resmin}
\end{eqnarray}
\begin{remark}
We define the discrete test space $V_m$ in such a way that it is as close as possible to the abstract $V$ space, to ensure stability, in a sense that the discrete inf-sup condition is satisfied. In our method it is possible to gain stability enriching the test space $V_m$ without changing the trial space $U_h$.
\end{remark}

Keeping in mind our definitions of bilinear form (\ref{eq:b}), right-hand-side (\ref{eq:l}) and the inner product (\ref{eq:inner}) our problem minimization of the residual is defined as:

Find $(r_m,u_h)_{V_m\times U_h}$ such as 
\begin{eqnarray}
(r_m,v_m)_{V_m} - 
\left(u_h,v_m\right) + \Delta t \left( \left(\beta_{\underline{i}} \frac {\partial u_h}{\partial x_{\underline{i}}}, v_m\right) + \left(\alpha\frac{\partial u_h}{\partial x_{\underline{i}}}, \frac{\partial v_m}{\partial x_{{\underline{i}}}} \right)\right) 
=  \nonumber  \\
\qquad - \Delta t\left(\left(\beta_{\underline{j}} \frac {\partial w}{\partial x_{\underline{j}}}, v_m\right) - \left(\alpha\frac{\partial w}{\partial x_{\underline{j}}}, \frac{\partial v_m}{\partial x_{\underline{j}}}\right)\right) +\left(w+f,v_m\right)\nonumber\\
 \quad \forall v_m\in V_m  \nonumber \\
\left(w_h,r_m\right) + \Delta t \left( \left(\beta_{\underline{i}} \frac {\partial w_h}{\partial x_{\underline{i}}}, r_m\right) + \left(\alpha\frac{\partial w_h}{\partial x_{\underline{i}}}, \frac{\partial r_m}{\partial x_{\underline{i}}}\right)\right) = 0 \quad \forall w_h \in U_h.
\label{eq:resmin}
\end{eqnarray}

We approximate the solution with tensor product of one dimensional B-splines basis functions of order $p$
\begin{equation}
u_h = \sum_{i,j} u_{i,j} B^x_{i;p}(x)B^y_{j;p}(y).
\end{equation}
We test with tensor product of one dimensional B-splines basis functions, where we enrich the order in the direction of the $x$ axis from $p$ to $q$ ($q\geq p$, and we enrich the space only in the direction of the alternating splitting)
\begin{equation}
v_m \leftarrow B^x_{i;q}(x)B^y_{j;p}(y).
\end{equation}
We approximate the residual with tensor product of one dimensional B-splines basis functions of order $p$
\begin{equation}
r_m = \sum_{s,t} r_{s,t} B^x_{s;q}(x)B^y_{t;p}(y),
\end{equation}
and we test again with tensor product of 1D B-spline basis functions of order $q$ and $p$, in the corresponding directions
\begin{equation}
w_h \leftarrow B^x_{k;p}(x)B^y_{l;p}(y).
\end{equation}
\begin{remark}
We perform the enrichment of the test space also in the alternating directions manner. In this way, when we solve the problem with derivatives along the $x$ direction, we enrich the test space by increasing the B-splines order in the $x$ direction, but we keep the B-splines order along $y$ constant (same as in the trial space).
By doing that, we preserve the Kronecker product structure of the matrix, to ensure that we can apply the alternating direction solver. 
\end{remark}

Let us focus on the particular entries of the left-hand-side matrix (\ref{eq:resmin}) discretized with tensor product B-spline basis functions. 
From the definition of the inner product in $V$
\begin{eqnarray}
\begin{aligned}
A_{k,l,s,t} =& \left( B^x_{s;q}(x)B^y_{t;p}(y) , B^x_{k;q}(x)B^y_{l;p}(y) \right)_{L_2(\Omega)} \nonumber \\
+& \left( B^x_{s;q-1}(x)B^y_{t;p}(y) , B^x_{k;q-1}(x)B^y_{l;p}(y) \right)_{L_2(\Omega)}.
\label{eq:A1}
\end{aligned}
\end{eqnarray}
From the definition of the $B$ operator
\begin{eqnarray}
\begin{aligned}
B_{k,l,i,j} =& \left( B^x_{i;p}(x)B^y_{j;p}(y) , B^x_{k;q}(x)B^y_{l;p}(y) \right)_{L_2(\Omega)}
\nonumber \\ +& \Delta t\left( \beta_x B^x_{i;p-1}(x)B^y_{j;p}(y) , B^x_{k;q}(x)B^y_{l;p}(y) \right)_{L_2(\Omega)}
\nonumber \\ +& \Delta t\left( \alpha B^x_{i;p-1}(x)B^y_{j;p}(y) , B^x_{k;q-1}(x)B^y_{l;p}(y) \right)_{L_2(\Omega)}.
\label{eq:B1}
\end{aligned}
\end{eqnarray}

We perform the direction splitting for all the terms:
\begin{eqnarray}
\begin{aligned}
A_{k,l,s,t} =& 
\left( B^x_{s;q}(x), B^x_{k;q}(x)\right)_{L_2(\Omega)}
\left( B^y_{t;p}(y) , B^y_{l;p}(y) \right)_{L_2(\Omega)}
 \nonumber \\
 + &
\left( B^x_{s;q-1}(x), B^x_{k;q-1}(x) \right)_{L_2(\Omega)}
\left( B^y_{t;p}(y) ,B^y_{l;p}(y) \right)_{L_2(\Omega)}
\label{eq:A1}\nonumber \\
= & 
\left(\left( B^x_{s;q}(x), B^x_{k;q}(x)\right)_{L_2(\Omega)}
+ \left( B^x_{s;q-1}(x), B^x_{k;q-1}(x) \right)_{L_2(\Omega)}\right) 
\nonumber \\
& \left( B^y_{t;p}(y) , B^y_{l;p}(y) \right)_{L_2(\Omega)},
\end{aligned}
\end{eqnarray}
\begin{eqnarray}
\begin{aligned}
B_{k,l,i,j} =& \left( B^x_{i;p}(x), B^x_{k;q}(x) \right)_{L_2(\Omega)}
\left(B^y_{j;p}(y), B^y_{l;p}(y) \right)_{L_2(\Omega)}
\nonumber\\ +& \Delta t
\left( \beta_x B^x_{i;p-1}(x), B^x_{k;q}(x)\right)_{L_2(\Omega)}
\left( B^y_{j;p}(y), B^y_{l;p}(y) \right)_{L_2(\Omega)}
\label{eq:B1}  \nonumber \\
+& \Delta t \alpha
\left( B^x_{i;p-1}(x), B^x_{k;q-1}(x)\right)_{L_2(\Omega)}
\left( B^y_{j;p}(y), B^y_{l;p}(y) \right)_{L_2(\Omega)}
\label{eq:B1} \nonumber \\
= &
\left(
\left( B^x_{i;p}(x), B^x_{k;q}(x) \right)_{L_2(\Omega)}
+ \Delta t\beta_x
\left( B^x_{i;p-1}(x), B^x_{k;q}(x)\right)_{L_2(\Omega)}
\label{eq:B1}  \right. \nonumber\\
+ &\left. \Delta t 
\left( \alpha B^x_{i;p-1}(x), B^x_{k;q-1}(x)\right)_{L_2(\Omega)}
\right) 
\left( B^y_{j;p}(y), B^y_{l;p}(y) \right)_{L_2(\Omega)}.
\end{aligned}
\end{eqnarray}

Our matrices are decomposed into two Kronecker product sub-matrices
\begin{equation}
A=A_x \otimes A_y; B=B_x \otimes B_y; B^{T}=B^{T}_x \otimes B^{T}_y; A_y=B_y
\end{equation}

\begin{equation}
\left(
\begin{array}{ll}
A & B \\
B^{T} & 0 \\
\end{array}
\right)
=
\left(
\begin{array}{ll}
A_x \otimes A_y & B_x \otimes B_y \\
B_x^{T}\otimes B_y^{T} & 0 \\
\end{array}
\right) = 
\end{equation}
\begin{equation}
=\left(
\begin{array}{ll}
A_x \otimes A_y & B_x \otimes A_y \\
B_x^{T}\otimes A_y & 0 \\
\end{array}
\right) = 
\left(
\begin{array}{ll}
A_x  & B_x \\
B_x^{T} & 0 \\
\end{array}
\right) \otimes A_y 
\end{equation}

Both matrices $\left(
\begin{array}{ll}
A_x & B_x \\
B_x^{T} & 0 \\
\end{array}
\right)$ and $A_y $ can be factorized in a linear ${\cal O}(N)$ computational cost.

\begin{figure}
\includegraphics[scale=0.4]{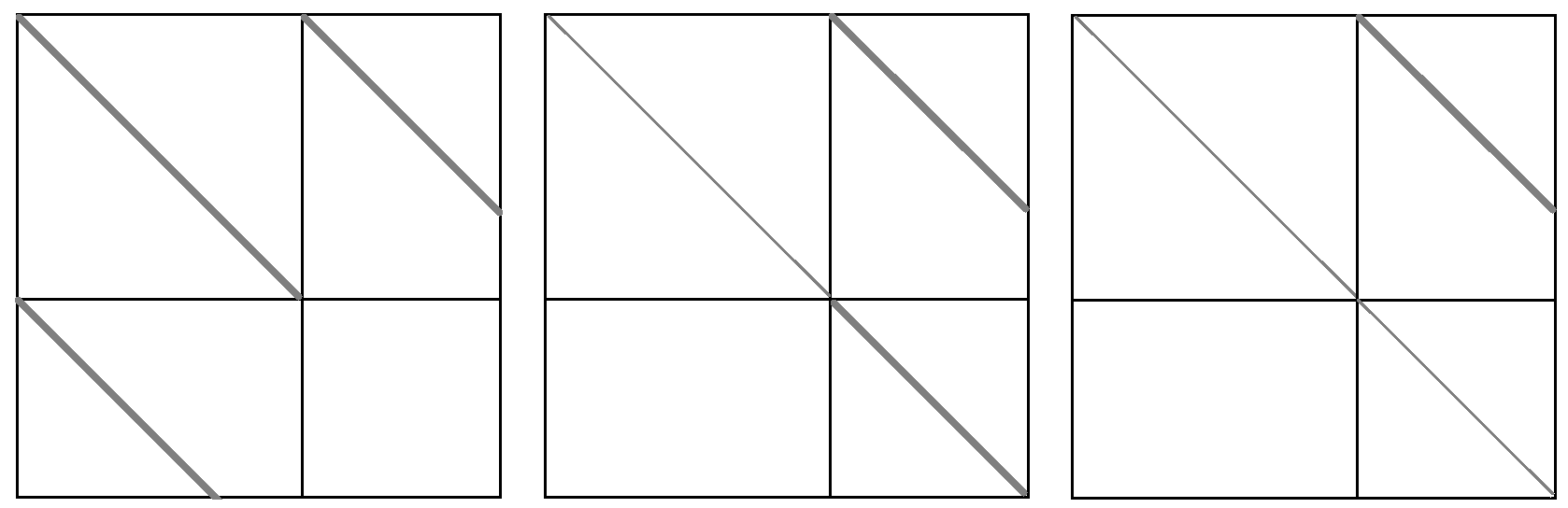}
\caption{Factorization of second block.}
\label{fig:block}
\end{figure}

The first block $\left(
\begin{array}{ll}
A_x & B_x \\
B_x^{T} & 0 \\
\end{array}
\right)$  consists of three multi-diagonal sub-matrices. The left-top block $A_x$ is of the size of $(N_x+r) \times (N_x+r)$, and it has a structure of the 1D Gramm matrix with B-splines of order $r$, so it has $2r+1$ diagonals.
The right-top block $B_x$ has a dimension of $(N_x+r) \times (N_x+p)$ and it is also a multi-diagonal matrix with $(r+p+1)$ diagonals. The left bottom part is just a transpose of the right-top one. This is illustrated in Figure \ref{fig:block}. The factorization uses the left-top block to zero the left-bottom block first. The computational cost of that is $(N_x+p)(r+p+1)={\cal }(N_x)$. Then, it makes the left-top block upper triangular. This is of the order of $(N_x+r)(2r+1)={\cal }(N_x)$. As the results of the first stage, in the right-bottom part we get the multi-diagonal sub-matrix with $(r+p+1)$ diagonals. We can continue the factorization of this right-bottom part, which costs $(N_x+p)(r+p+1)={\cal }(N_x)$. At the end, we run the global backward substitution, which costs 
$(N_x+p)(r+p+1)+(N_x+r)(r+p+1+2r+1)={\cal O}(N_x)$.

Identical consideration is taken for the second sub-step, resulting in linear computational cost alternating directions solver for implicit time integration scheme.

\section{Numerical results}

\subsection{Verification of the time-dependent simulations}

In order to verify the time-integration schemes and the direction splitting solver, we construct a time-depedent advection-diffusion problem with manufactured solution. Namely, we consider
two-dimensional problem
\begin{equation}
\frac{d u}{dt} -\nabla \cdot \left(\alpha \nabla u\right) + \beta \cdot \nabla u = f, \nonumber
\end{equation}
with $\alpha=10^{-2}$, $\beta=(1,0)$, with zero Dirichlet boundary conditions solved on a square $[0,1]^2$ domain.
We setup the forcing function $f(x,y,t)$ in such a way that it delivers the manufactured solution of the form $u(x,y,t)=\sin(\pi x)\sin(\pi y)\sin(\pi t)$ on a time interval $[0,2]$.

We solve the problem with residual minimization method on $10 \times 10$ mesh with 100 time steps $\Delta t=0.01$, using the Kronecker product solver with linear computational cost.

We compute the error between the exact solution $u_{exact}$ and the numerical solution $u_h$. We present the comparisons with different time step size $\Delta t$. We compute
We compute relative $L^2$ norm $\|u_{\textrm{exact}}(t) - u_{\textrm{h}}(t)\|_{L^2} / \|u_{\textrm{exact}}(t)\|_{L^2} * 100 \%$, and relative $H^1$ norm
$\|u_{exact}(t) - u_{\textrm{h}}(t)\|_{H^1} / \|u_{\textrm{exact}}(t)\|_{H^1} * 100 \%$.

\begin{figure}
\includegraphics[scale=0.65]{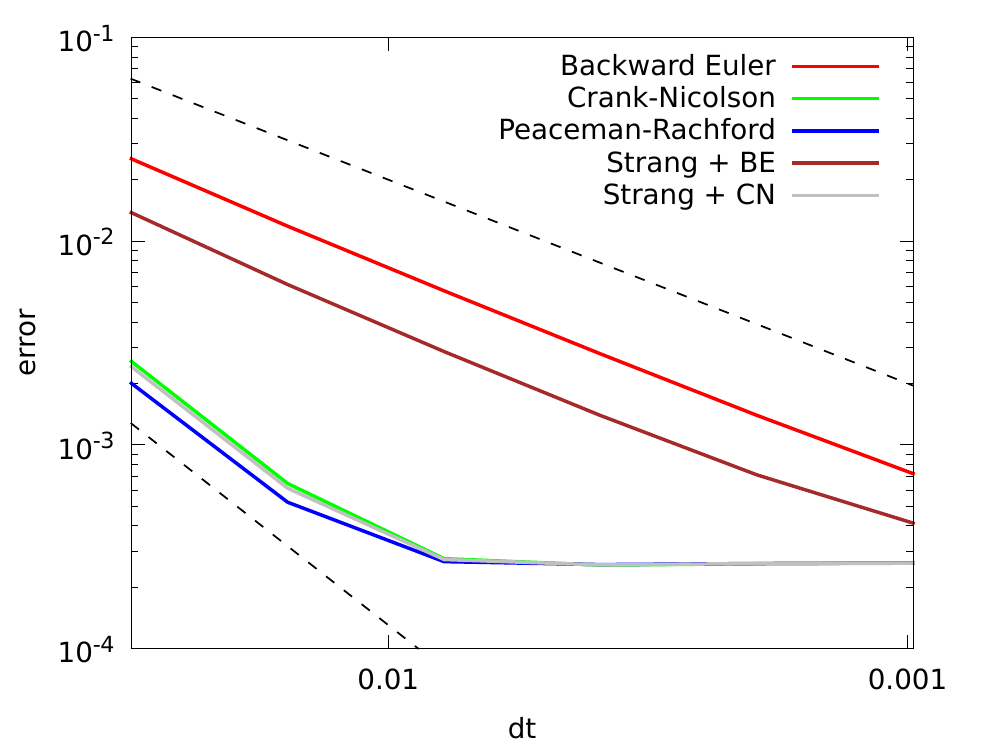}
\includegraphics[scale=0.65]{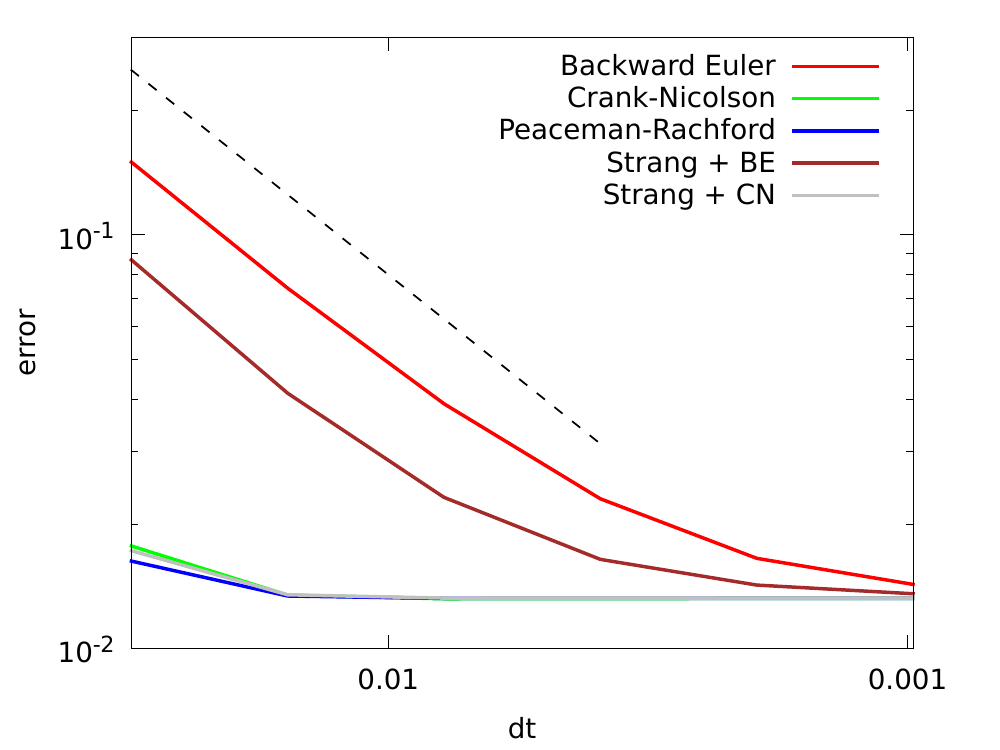}
\caption{Convergence in $L^2$ and $H^1$ norms for different time integration scheme on $8\times 8$ mesh.}
\label{fig:L2_8}
\end{figure}

\begin{figure}
\includegraphics[scale=0.65]{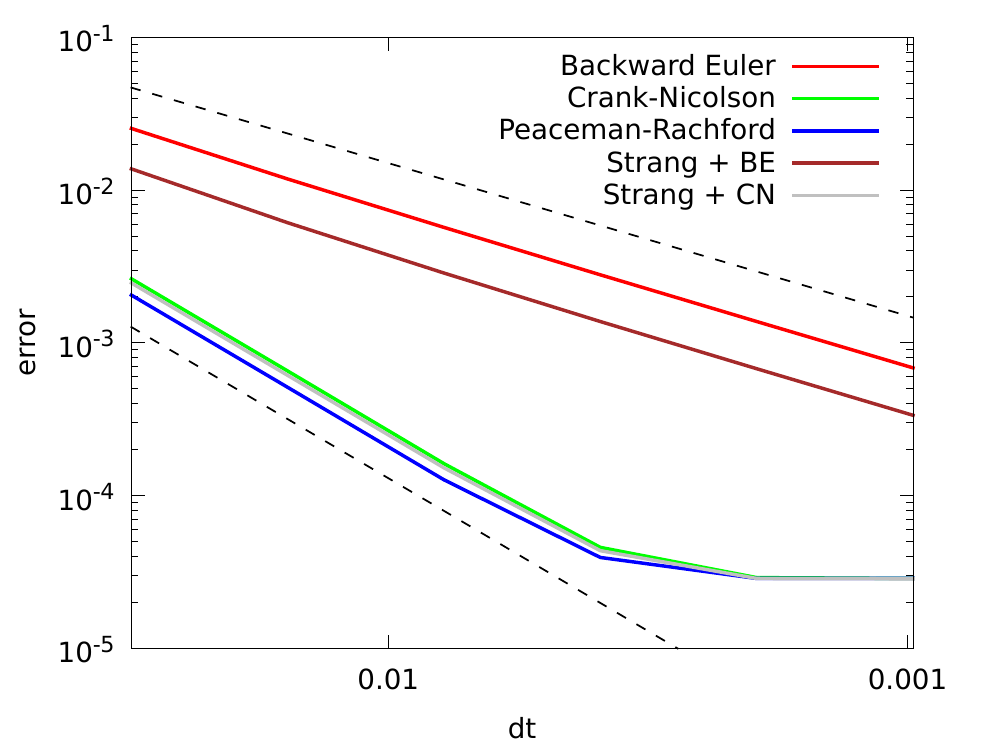}
\includegraphics[scale=0.65]{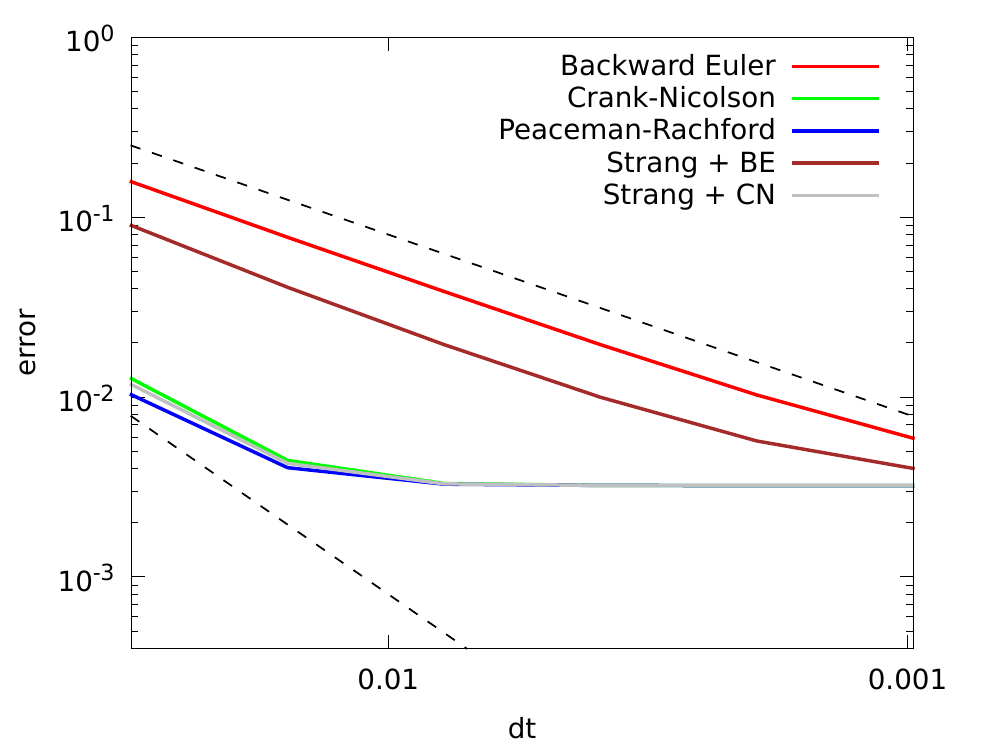}
\caption{Convergence in $L^2$ and $H^1$ norms for different time integration scheme on $16\times 16$ mesh.}
\label{fig:L2_16}
\end{figure}

\begin{figure}
\includegraphics[scale=0.65]{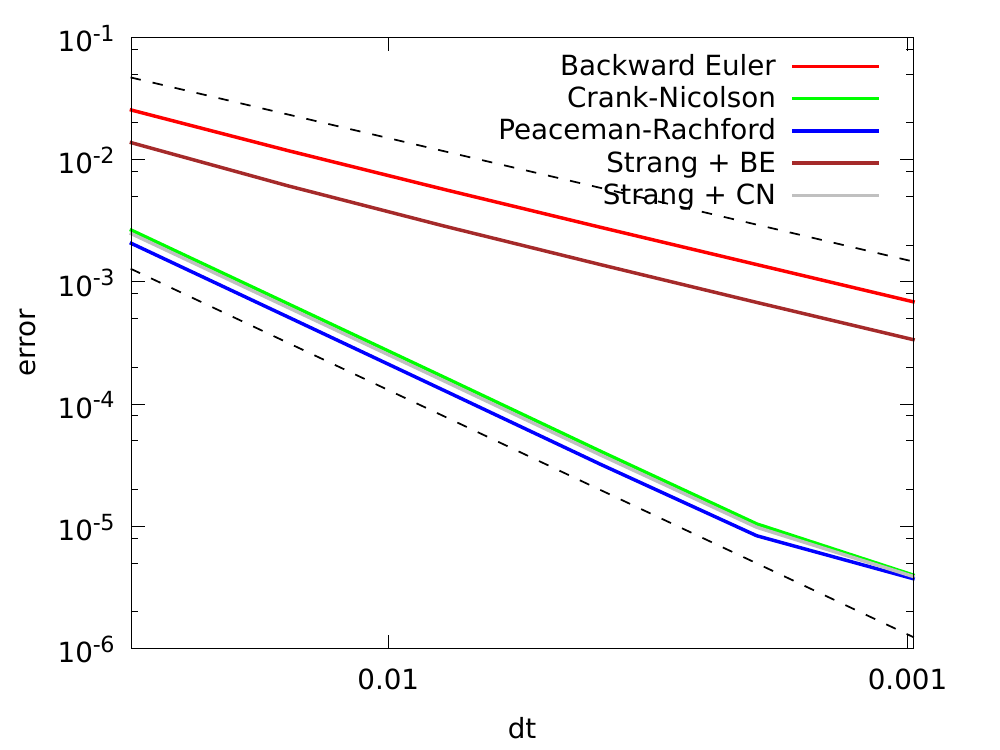}
\includegraphics[scale=0.65]{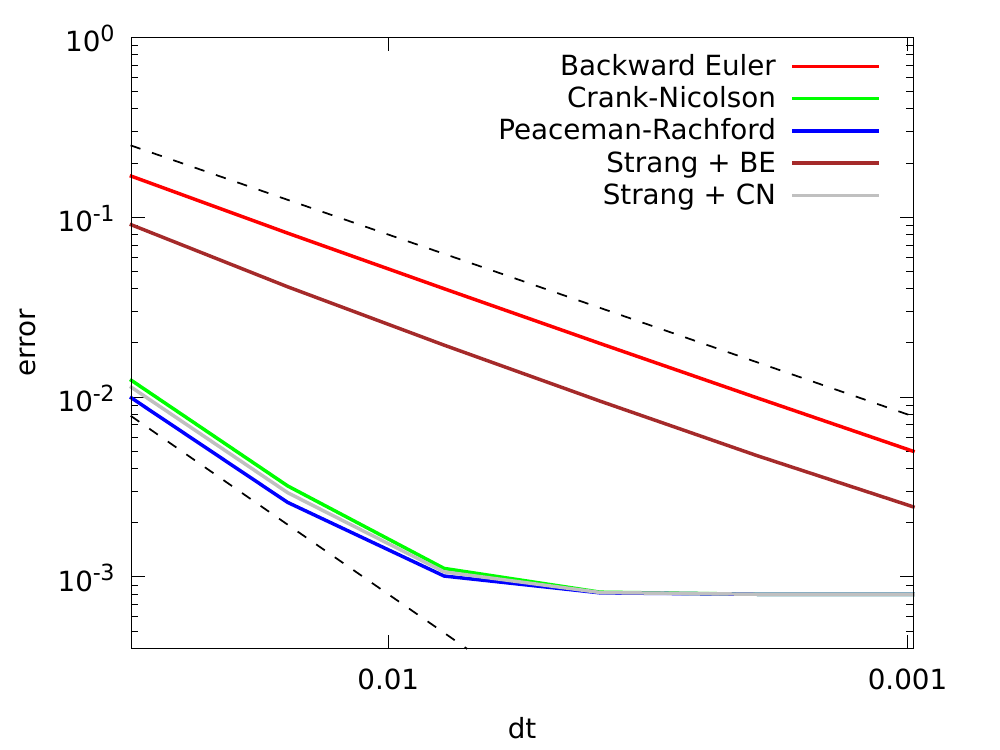}
\caption{Convergence in $L^2$ and $H^1$ norm for different time integration scheme on $32\times 32$ mesh.}
\label{fig:L2_32}
\end{figure}

The $L^2$ and $H^1$ errors for different mesh dimensions, for trial (2,1) test (2,0) are presented in Figures \ref{fig:L2_8}-\ref{fig:L2_32}. We can read the following conclusions from these experiments:
\begin{itemize}
\item The backward Euler scheme itself as well as Strang method with the backward Euler schere are of the first order.
\item The Peacemen-Rachford scheme is of the second order.
\item The Crank-Nicolson scheme is of the second order.
\item The Strang method with Crank-Nicolson scheme is the second order scheme as well.
\item The maximum time-integration scheme accuracy depends on the spatial accuracy, and for the mesh size $32\times 32$ it is of the order of $10^{-5}$ for the second order time integration schemes.
\end{itemize}


\begin{figure}
\includegraphics[scale=0.65]{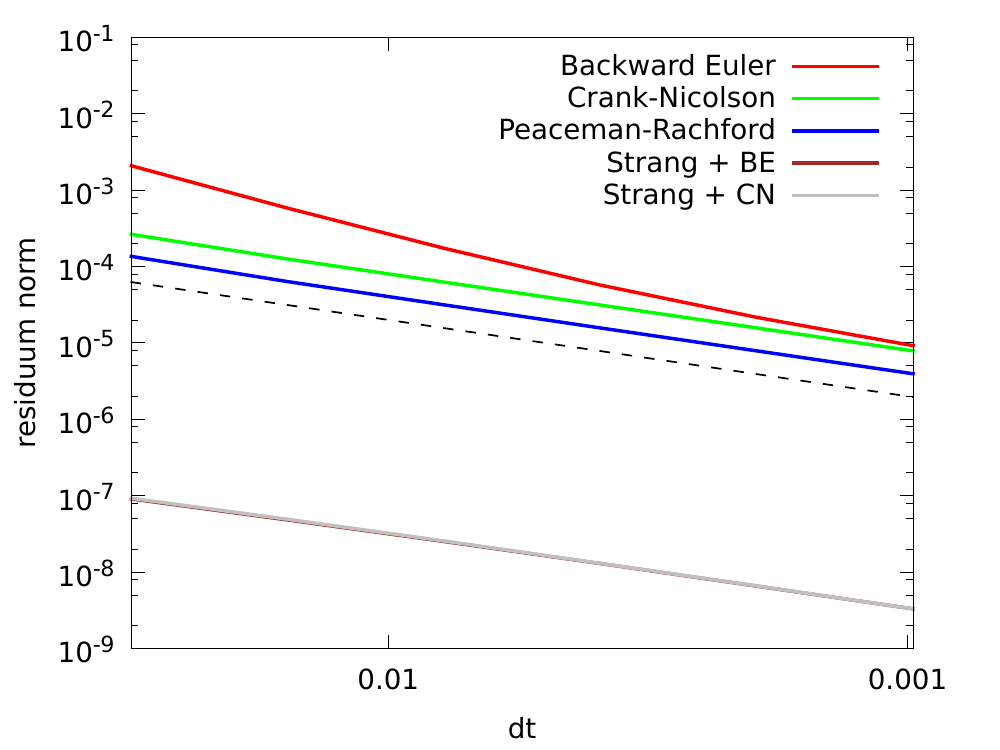}
\includegraphics[scale=0.65]{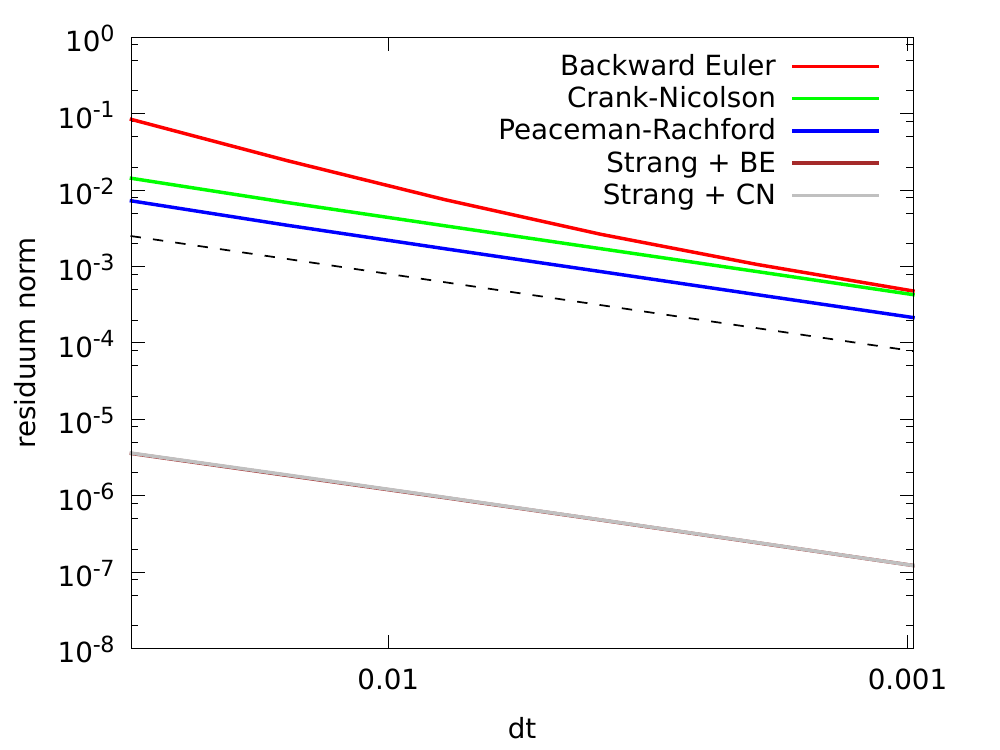}
\caption{$L^2$ and $H^1$ norms of the residual for different time integration scheme on $8\times 8$ mesh.}
\label{fig:r_L2_8}
\end{figure}

\begin{figure}
\includegraphics[scale=0.65]{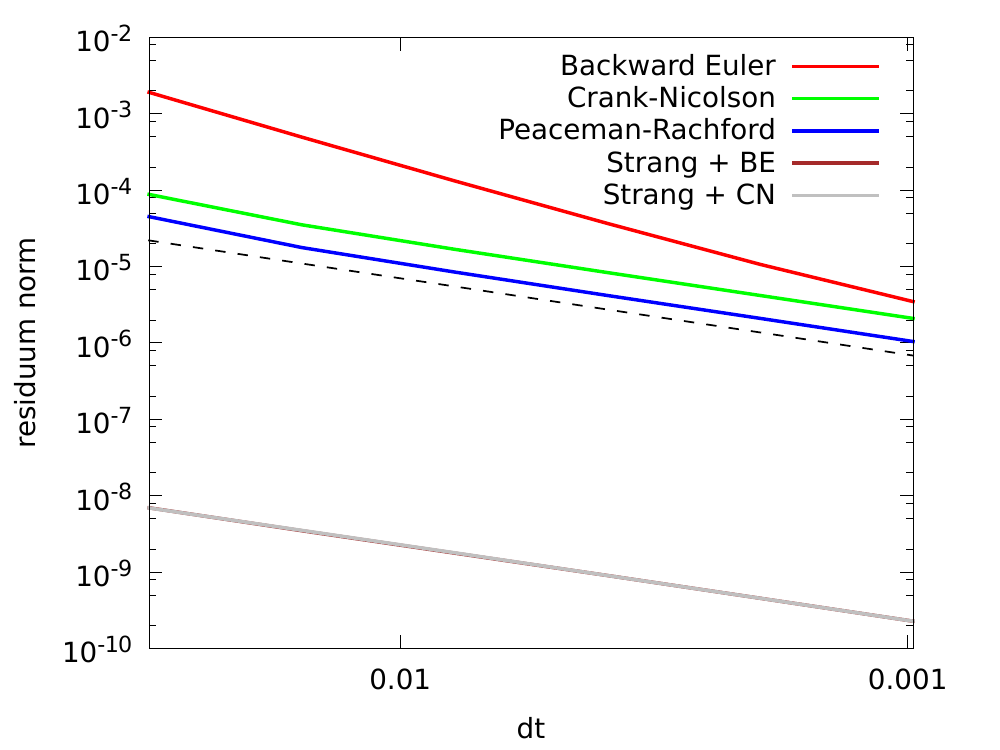}
\includegraphics[scale=0.65]{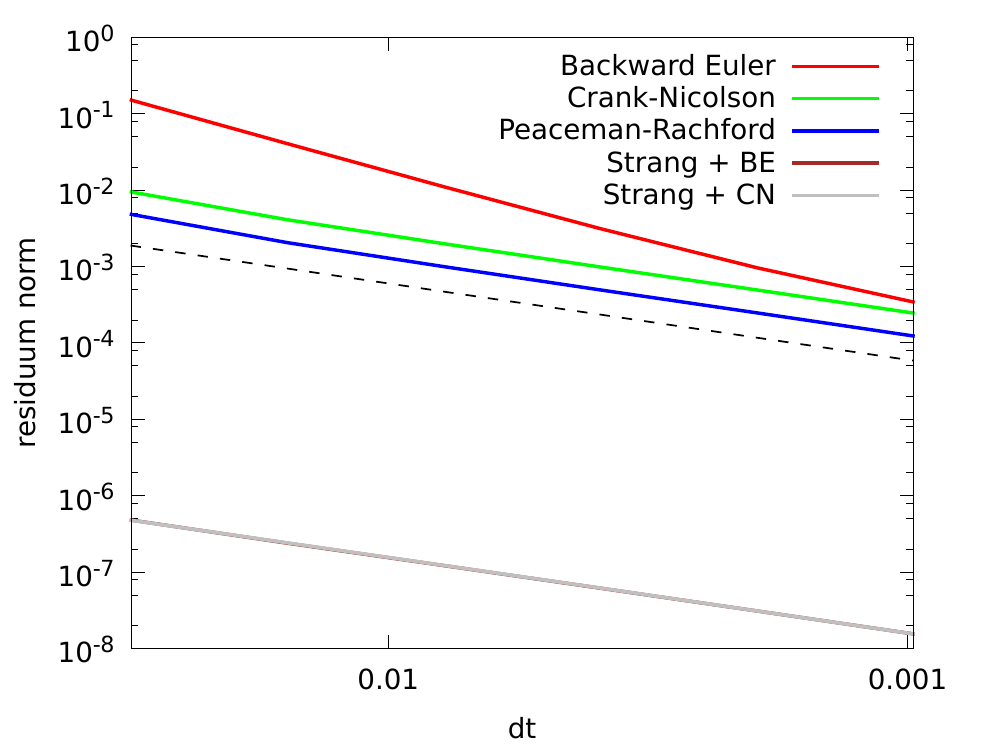}
\caption{$L^2$ and $H^1$ norms of the residual for different time integration scheme on $16\times 16$ mesh.}
\label{fig:r_L2_16}
\end{figure}

\begin{figure}
\includegraphics[scale=0.65]{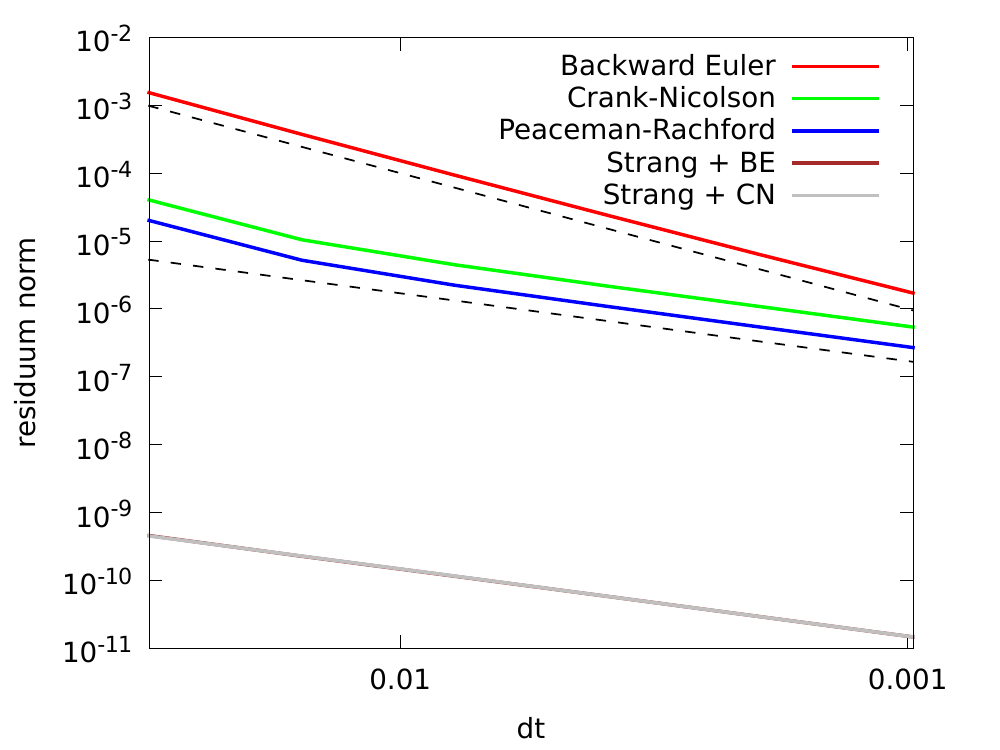}
\includegraphics[scale=0.65]{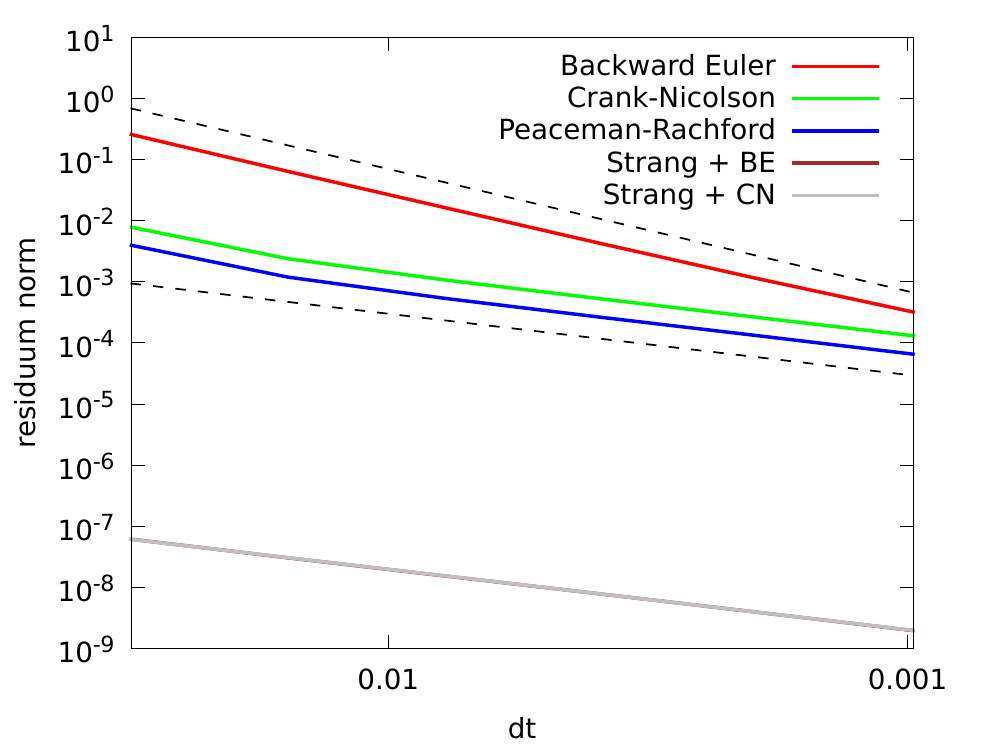}
\caption{$L^2$ and $H^1$ norm of the residual for different time integration scheme on $32\times 32$ mesh.}
\label{fig:r_L2_32}
\end{figure}

The residual minimization method allows using the norm of the residual computed during the solution process as the error estimation for possible adaptivity. We present in Figures \ref{fig:r_L2_8}-\ref{fig:r_L2_32} the plots of the residual $L^2$ and $H^1$  norms, for all the proposed time integration schemes. We can see that these residual norms behave in a different way than the already presented errors computed by using the known exact solution. Namely, the order of the method can be checked by looking at the exact errors, and not by looking at the residual errors. However, these residual errors give a good estimate for the adaptive procedure. Moreover, the two versions of the Strang splitting method, as measured with the residual, give identical plots.


\begin{figure}[h]
\includegraphics[scale=0.2]{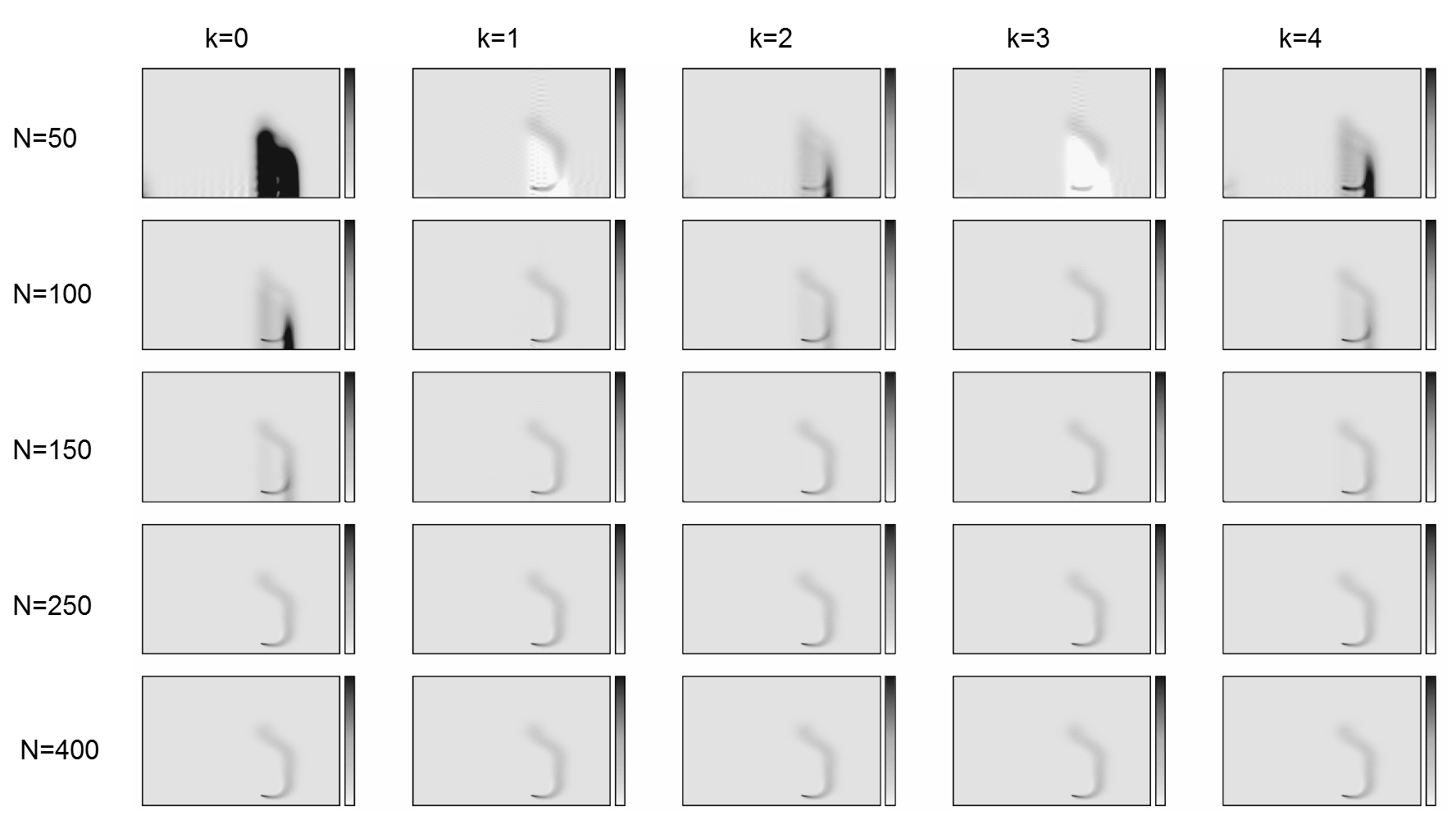}
\caption{Stabilization of 2D simulation with quadratic B-splines for trial spaces with iGRM method. Different rows correspond to different mesh sizes, varying from $N=50, 100, 150, 250, 400$. Different columns correspond to different order of enrichment of test B-splines of order $p+k$ and continuity $p+k-1$.}
\label{fig:2D}
\end{figure}

\subsection{Propagation of the pollutant from a chimney}
In this section we present numerical results for the advection-diffusion equation in two-dimensions. 
We use the pollution model based on \cite{Pollution}, with the following simplifying assumptions: we consider only one component of the pollution vector field, we neglect the chemical interactions between different components, and we assume cube shape of the computational domain. These assumptions have not decreased the complexity of the stabilization issues related to the numerical simulations. 
Hence we propose to test our method on this simple but reliable model, where the Galerkin discretization is known to be unstable.

We utilize tensor product of 1D B-splines along the $x$, and $y$ and we perform the direction splitting with three intermediate time steps, resulting in the three linear systems having Kronecker product structure. In our equation we have now $\epsilon=(50+sin(\frac{x\pi}{5000}),50+\frac{y}{5000})$ in 2D. The computational domain is a cube with dimensions of $5000\times 5000$ meters in 2D

The wind $\beta$ is given by
\begin{equation}
\beta(x,y)=(\cos a(t),\sin a(t)),
\end{equation}
where 
\begin{equation}
a(t)=\pi/3 (\sin(s) + 0.5 \sin(2.3 s)) + 3/8 \pi,
\end{equation}
where $s=t/150$.

The right-hand-side represents a chimney, a source of pollution given by
\begin{equation}
f(p) = (r - 1)^2 (r + 1)^2 ,
\end{equation}
where $r = \min(1, (|p - p_0| / 25)^2)$, $p$ represents the distance from the chimney, and $p_0$ is the chimney location.
The initial state was the ambient concentration of the pollution, of the order of $10^{-6}$ in the entire domain.

We solve the problem with residual minimization method on two-dimensional meshes with different dimensions, using the Kronecker product solver with linear computational cost.

We utilize trial space (2,1) or (3,2) and we experiment with (2,1), (3,2), (4,3), (5,4), (6,5) or (2,0), (3,0), (3,1) test spaces, as well as we try trial space (3,2) and we experiment with (3,1) or (3,1) test spaces.

\begin{figure}[h]
\includegraphics[scale=0.2]{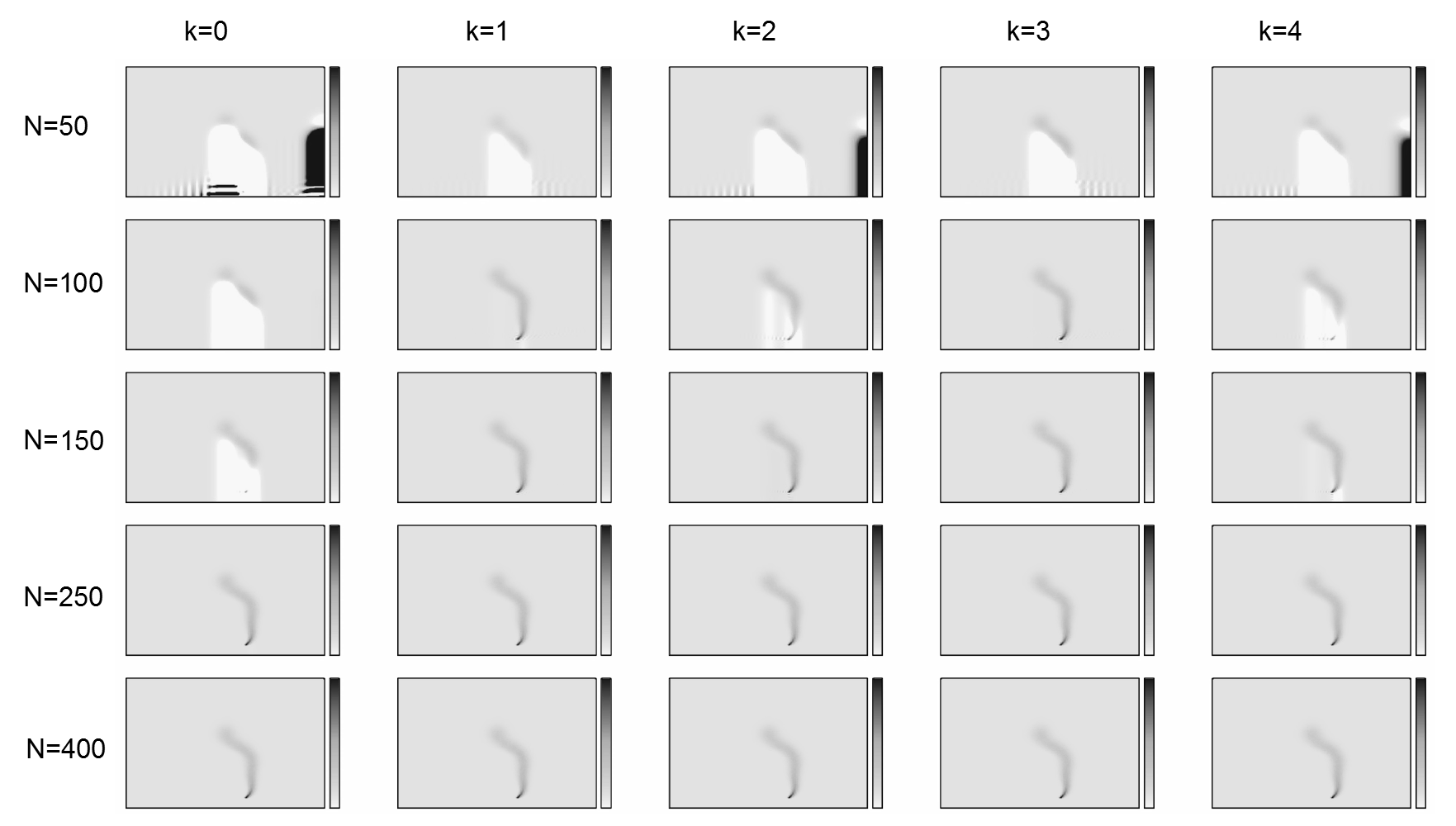}
\caption{Stabilization of 2D simulation with cubic B-splines for trial spaces with iGRM method. Different rows correspond to different mesh sizes, varying from $N=50, 100, 150, 250, 400$. Different columns correspond to different order of enrichment of test B-splines of order $p+k$ and continuity $p+k-1$.}
\label{fig:2Db}
\end{figure}

\begin{figure}[h]
\includegraphics[scale=0.4]{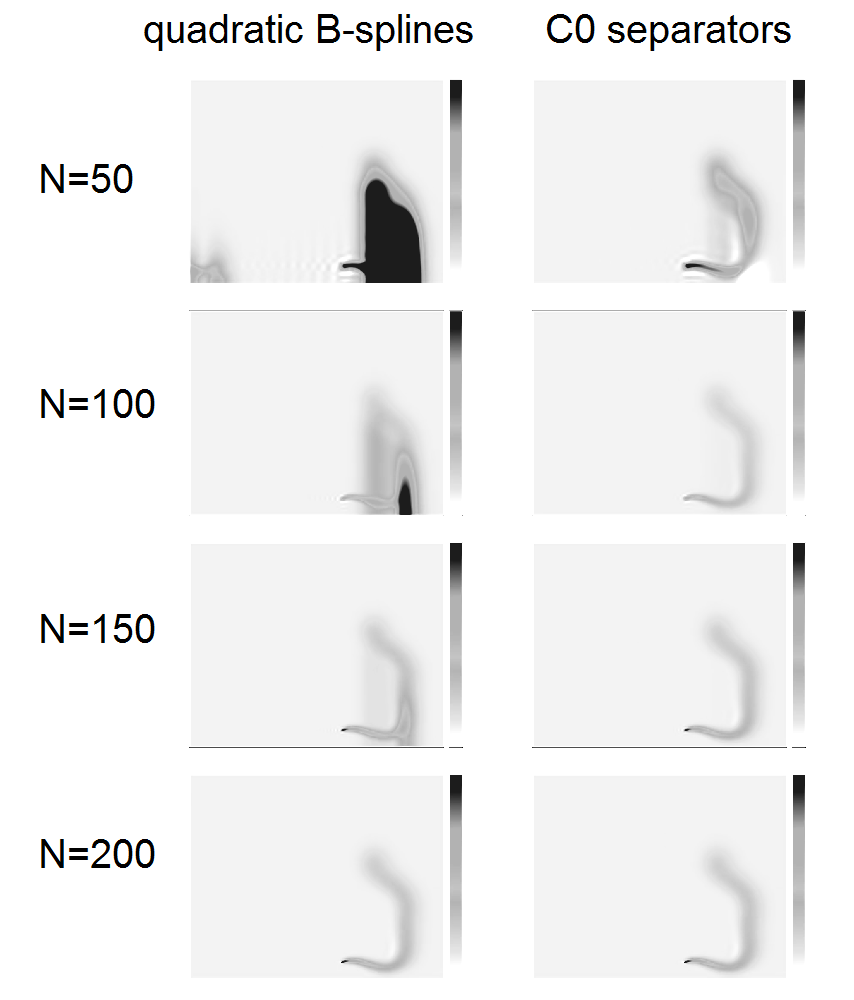}
\caption{Stabilization of 2D simulation with quadratic B-splines for trial spaces with iGRM method. Different rows correspond to different mesh sizes, varying from $N=50, 100, 150, 250, 400$. The first column corresponds to the Galerkin method with quadratic B-splines. The second column corresponds to the test space enriched by introduction of $C^0$ separators between all the elements, which is equivalent to the Lagrange polynomials.}
\label{fig:2Dc}
\end{figure}

\begin{figure}[h]
\includegraphics[scale=0.2]{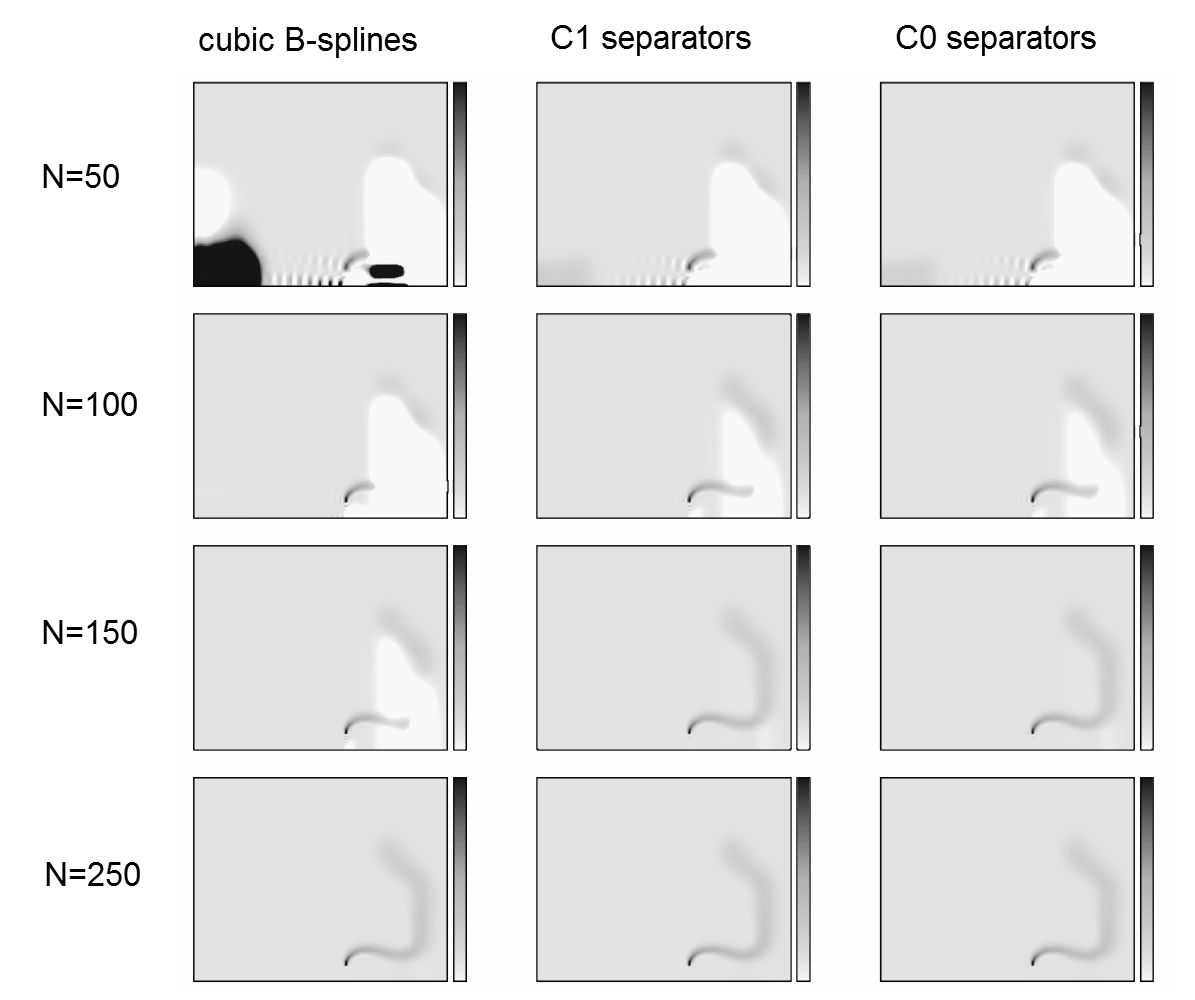}
\caption{Stabilization of 2D simulation with cubic B-splines for trial spaces with iGRM method. Different rows correspond to different mesh sizes, varying from $N=50, 100, 150, 250, 400$. The first column corresponds to the Galerkin method. The second column corresponds to single repetition of knot points between all the elements (reduction of continuity down to $C^1$ between the elements).
The third columns corresponds to the test space enriched by introduction of $C^0$ separators between all the elements.}
\label{fig:2Dd}
\end{figure}

\begin{figure}[h]
\includegraphics[scale=0.2]{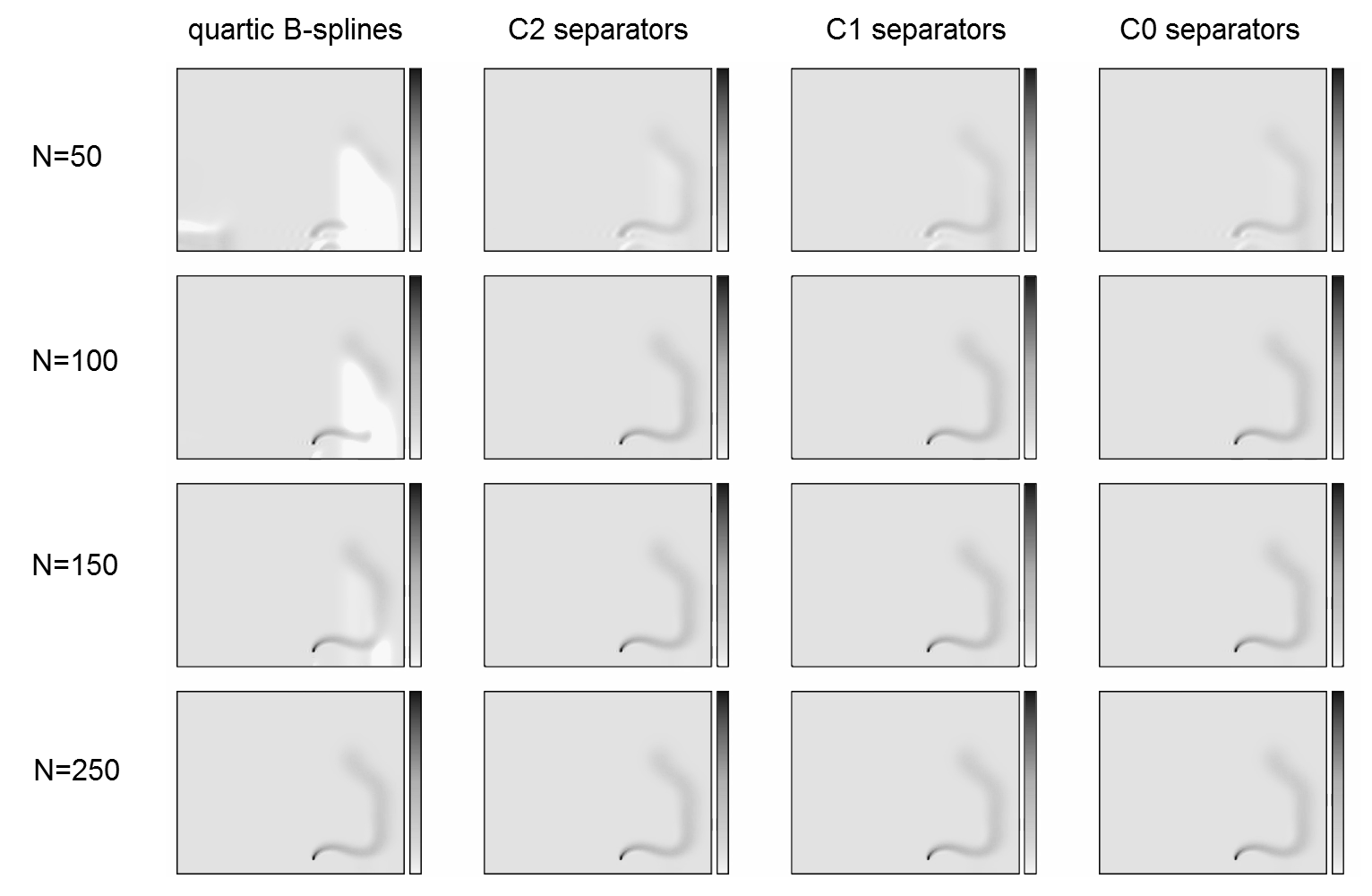}
\caption{Stabilization of 2D simulation with quartic B-splines for trial spaces with iGRM method. Different rows correspond to different mesh sizes, varying from $N=50, 100, 150, 250, 400$. The first column corresponds to the Galerkin method. The second column corresponds to single repetition of knot points between all the elements (reduction of continuity down to $C^2$ between the elements).
The second column corresponds to twofold repetition of knot points between all the elements (reduction of continuity down to $C^1$ between the elements).
The fourth columns corresponds to the test space enriched by introduction of $C^0$ separators between all the elements.}
\label{fig:2De}
\end{figure}

The 2D numerical results are presented in Figures \ref{fig:2D}-\ref{fig:2De}.
The Figures differ by the way we define the trial space, and by the way we perform the enrichment of the test space.

In Figures \ref{fig:2D}-\ref{fig:2Db} we used quadratic and cubic B-splines for the trial space. We performed 100-time steps of the implicit stable method with linear ${\cal O}(N)$ computational cost solver. We present snapshots of the numerical results for different mesh dimensions (in rows) and different B-spline orders of the enriched spaces (in columns). We can read that small mesh size ($N=100$) and $k=1$ (enrichment of the B-spline order of $1$) gives stable simulations, comparable to the results of the standard Galerkin method on $N=250$ mesh. 

In Figures \ref{fig:2Dc}-\ref{fig:2De} we use quadratic, cubic, and quartic B-splines for trial space. We present snapshots of the numerical results for different mesh dimensions (in rows). In columns, we repeat knot points between elements, so we reduce the continuity between the elements, down to $C^0$. Note that reducing the continuity actually \emph{enriches} the test space, since it is equivalent to introduction of a new basis functions.

\subsection{Circular wind problem}
Finally, we consider the circular wind problem, where the advection coefficients do not have the Kronecker product structure. In this case, we cannot utilize the linear computational cost direction splitting solver, but we still can call e.g. direct solver MUMPS \cite{MUMPS1,MUMPS2,MUMPS3} to factorize the residual minimization problem, or call a preconditioned conjugate gradients iterative solver. 

\begin{figure}
\includegraphics[scale=0.4]{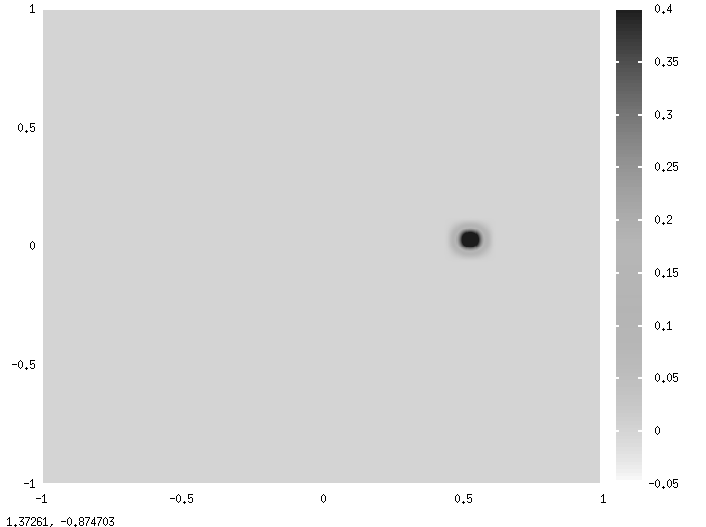}
\caption{Initial state for the simulation of the time-dependent advection-diffusion problem with circular wind and Pecklet number $Pe=1,000,000$.}
\label{fig:initial}
\end{figure}

Namely, we consider time-dependent advection diffusion equations 
\begin{equation}
\frac{\partial u}{\partial t}+\beta_x \frac{\partial u}{\partial x}+\beta_y \frac{\partial u}{\partial y}-\alpha \left(\frac{\partial^2 u}{\partial x^2}+
\frac{\partial^2 u}{\partial y^2}\right)=0
\label{eq:Problem2}
\end{equation}
with $\alpha=10^{-6}$ and the right-hand side representing the circular wind
\begin{equation}
\beta(x,y)=(y,-x)
\end{equation}
and the initial state presented in Figure \ref{fig:initial}.

\begin{figure}
\centering
\begin{subfigure}[b]{0.45\textwidth}
\includegraphics[scale=0.35]{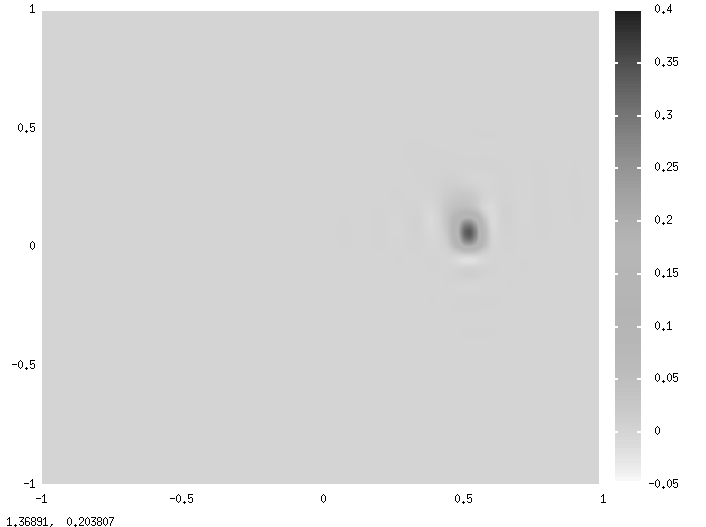}
\caption{t=1}
\end{subfigure}
\begin{subfigure}[b]{0.45\textwidth}
\includegraphics[scale=0.35]{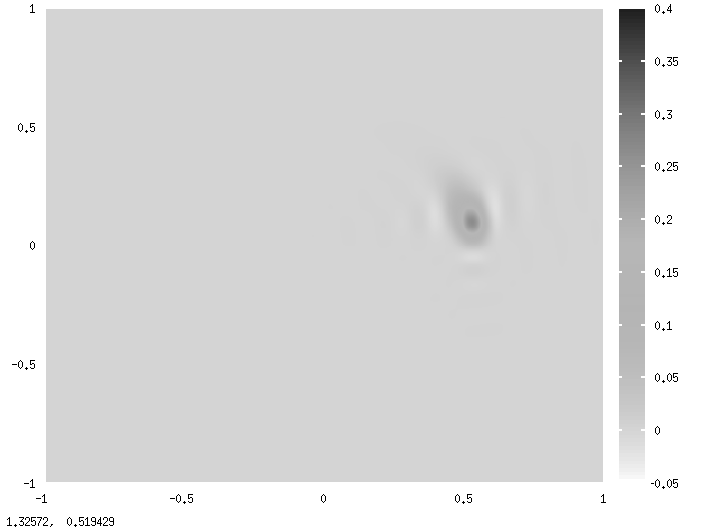}
\caption{t=2}
\end{subfigure}
\begin{subfigure}[b]{0.45\textwidth}
\includegraphics[scale=0.35]{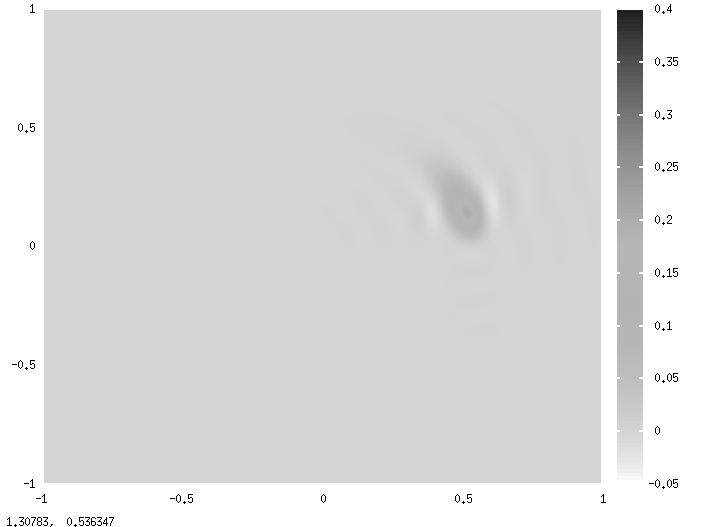}
\caption{t=3}
\end{subfigure}
\begin{subfigure}[b]{0.45\textwidth}
\includegraphics[scale=0.35]{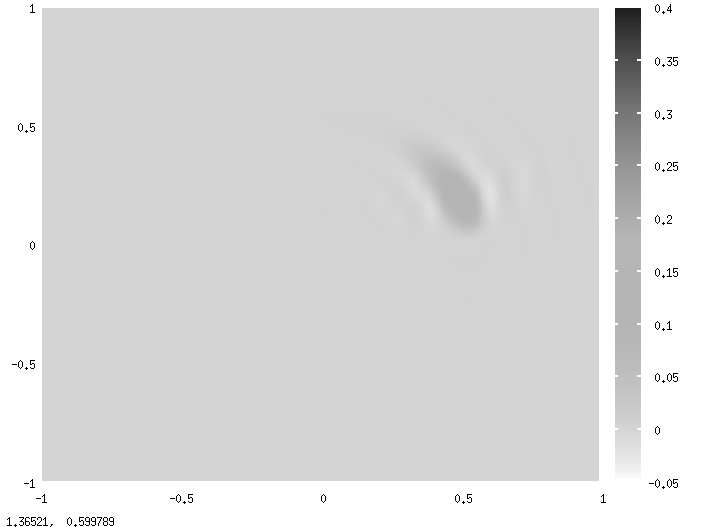}
\caption{t=4}
\end{subfigure}
\begin{subfigure}[b]{0.45\textwidth}
\includegraphics[scale=0.35]{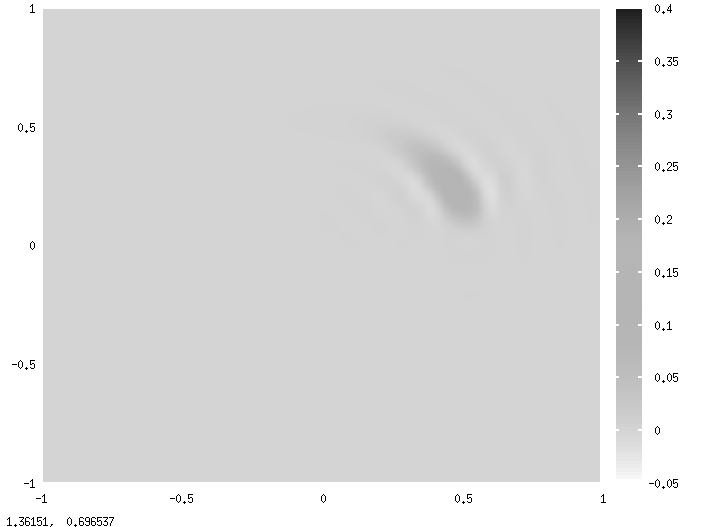}
\caption{t=5}
\end{subfigure}
\begin{subfigure}[b]{0.45\textwidth}
\includegraphics[scale=0.35]{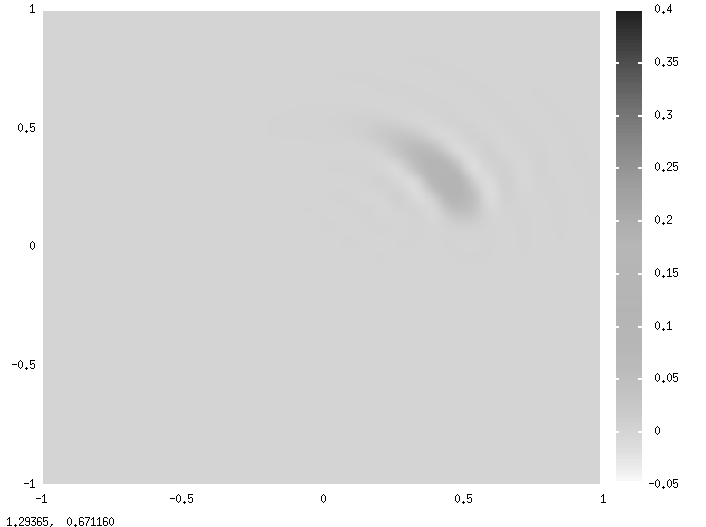}
\caption{t=6}
\end{subfigure}
\begin{subfigure}[b]{0.45\textwidth}
\includegraphics[scale=0.35]{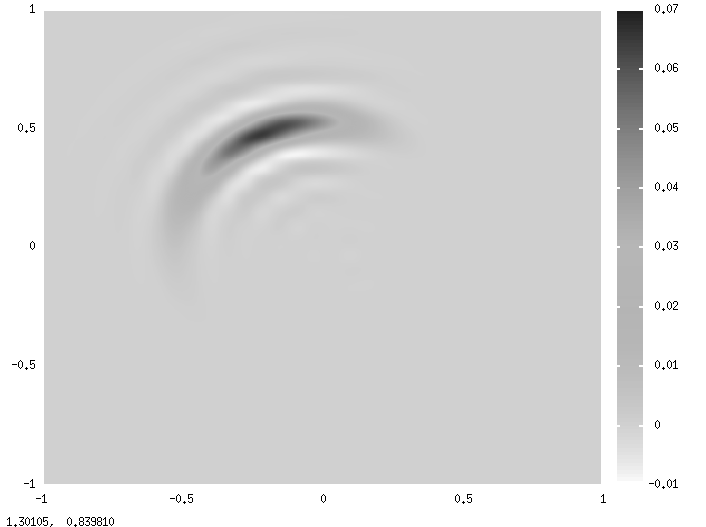}
\caption{t=20}
\end{subfigure}
\begin{subfigure}[b]{0.45\textwidth}
\includegraphics[scale=0.18]{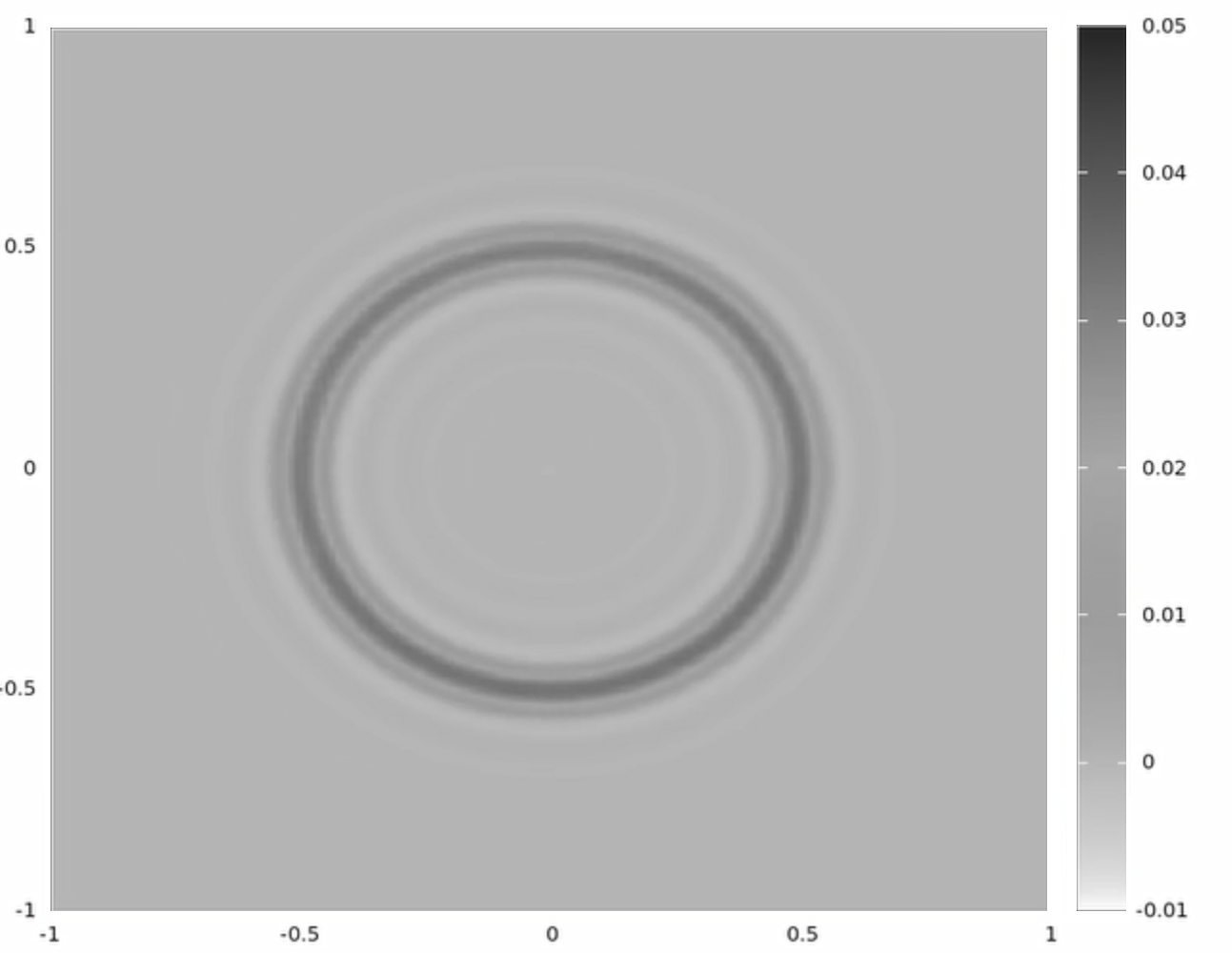}
\caption{t=1000}
\end{subfigure}
\caption{Snapshoots from the simulation of the time-dependent advection-diffusion problems with circular wind and Pecklet number $Pe=1,000,000$.}
\label{fig:circular_wind}
\end{figure}

We formulate the residual minimization system according to (\ref{eq:resmin}). 
To gain stability we utilize computational mesh with $64\times 64$ elements, with trial$=(4,3)$ (quartic B-splines with $C^3$ continuity) and test$=(5,0)$ (quintic B-splines with $C^0$ continuity).
We solve the iGRM system using the MUMPS solver at every time step. The snapshots from the simulations for $\alpha=10^{-6}$ are presented in Figures \ref{fig:circular_wind}.

\begin{figure}
\centering
\includegraphics[scale=0.25]{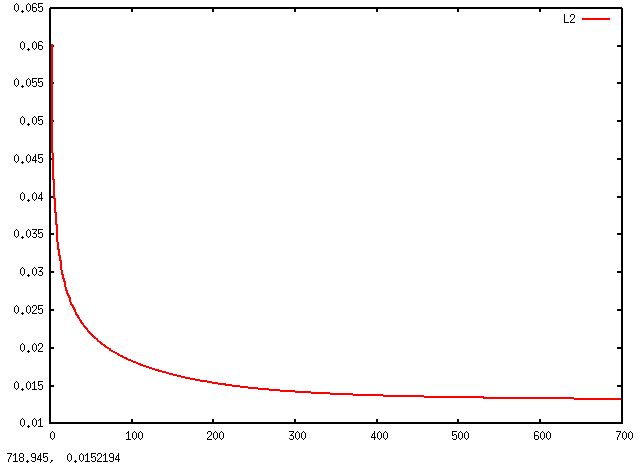}
\includegraphics[scale=0.25]{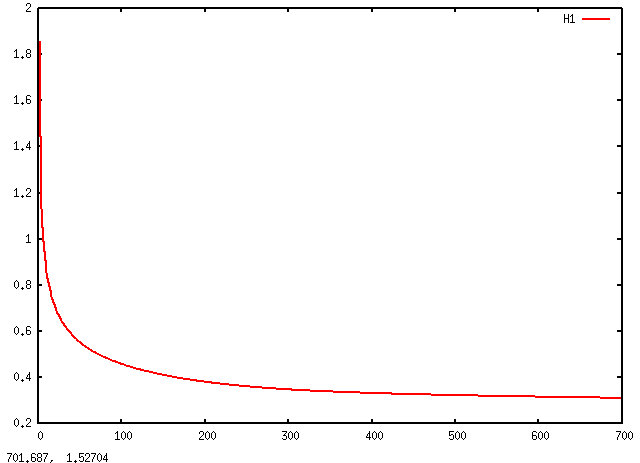}
\caption{$L^2$ and $H^1$ norms of the solution for the simulation of the time-dependent advection-diffusion problem with circular wind and Pecklet number $Pe=1,000,000$.}
\label{fig:L2time}
\end{figure}

The convergence of the $L^2$ and $H^1$  norms over the time are illustrated in Figure \ref{fig:L2time}.

Since the computational cost is no longer linear ${\cal O}(N)$, we provide in Table 1 the execution times for MUMPS solver called for different configurations of mesh dimensions and trial and test space.
\begin{table*}[htp]
\begin{center}
\begin{tabular}{|c|c|c|c|}
\hline 
space       & n        & dofs    & solver [ms] \\
\hline
trial = (2,1), test = (3,0)        &8        &725          &5 \\
trial = (2,1), test = (3,0)       &16       &2725         &27 \\ 
trial = (2,1), test = (3,0)       &32       &10565       &259 \\
trial = (2,1), test = (3,0)       &64       &41605       &1393 \\
\hline
trial = (3,2), test = (4,0)        &8        &1210         &14 \\ 
trial = (3,2), test = (4,0)       &16       &4586         &100 \\
trial = (3,2), test = (4,0)       &32       &17866       &928 \\
trial = (3,2), test = (4,0)       &64       &70538       &5725 \\
\hline
trial = (4,3), test = (5,0)       &8        &1825         &51 \\
trial = (4,3), test = (5,0)       &16      &6961         &618 \\
trial = (4,3), test = (5,0)       &32      &27217       &2722 \\
trial = (4,3), test = (5,0)       &64      &107665     &9116 \\
\hline
\end{tabular}
\caption{Trial and test spaces, mesh dimensions, total number of degrees of and execution time in [ms] as reported by MUMPS solver.}
\end{center}
\label{tab1}
\end{table*}

\section{Conclusions}
\label{sec:conc}

We introduced a stabilized isogeometric finite element method for the implicit problems that result in a Kronecker product structure of the computational problem. This, in turn, gives the linear computational cost of the implicit solver, and the smooth approximation of the solution resulting from the B-spline basis functions. We called our method isogeometric residual minimization (iGRM). The method has been verified on the three advection-diffusion problems, including the model ``membrane'' problem, the problem modeling the propagation of the pollution from a chimney, and the circular wind problem. Our future work will involve the development of this method for other equations, as well as the development of the parallel software dedicated to the simulations of different non-stationary problems with the iGRM method. We are going also to develop the mathematical fundations on the error analysis, and the convergence of the method. 

\subsection*{Acknowledgments}

This work is partially supported by National Science Centre, Poland grant no. 2017/26/M/ST1/00281.

\end{document}